# Šafarevič's Theorem on Solvable Groups as Galois Groups

by Alexander Schmidt and Kay Wingberg

The aim of this article is to give a complete proof of the following famous theorem of *I. R. Šafarevič*:

**Theorem 1** *Let $k$ be a global field and let $G$ be a finite solvable group. Then there exists a finite Galois extension $K|k$ with Galois group $G(K|k) \cong G$.*

Šafarevič proved this result in 1954. The intricate proof uses arguments which are spread over several papers and unfortunately, it contains a mistake relative to the prime 2. [1]

In the course of the years, several mathematicians tried to simplify this highly complicated proof using the new development in number theory since then. Nevertheless, so far no complete proof has been published.

In the following we want to present an up to date proof of theorem 1 which also corrects the mistake in the original article. Of course, we will use the original ideas of Šafarevič, in particular the remarkable technique of *shrinking* obstructions, which is highly instructive. The authors do not know any other argument in number theory, which utilizes a similar technique.

1. Introduction

The first attempt to prove a statement like theorem 1 was made in 1936 independently by *A. Scholz* and *H. Reichardt*, see [10]. Their approach applies only to groups of order prime to $\#\mu(k)$, the order of the group of roots of unity contained in the global field $k$. In particular, it can be used to prove that every (solvable) group of odd order is a Galois group over $\mathbb{Q}$. The reader can find a proof of this statement for nilpotent groups in the spirit of Scholz and Reichardt in Serre's book [16].

As shown by *J. Neukirch*, the method of Scholz and Reichardt can even be strengthened in order to construct global solvable extensions which realize finitely many given local extensions:

---

[1]Šafarevič in his collected papers explains shortly how to correct this and we will follow his suggestion here.

**Theorem (Neukirch [7])** *Let $k$ be an algebraic number field and let $S$ be a finite set of primes of $k$. Let $G$ be a pro-solvable group of finite exponent prime to $\#\mu(k)$ and for $\mathfrak{p} \in S$ let $K_\mathfrak{p}|k_\mathfrak{p}$ be Galois extensions whose Galois groups $G(K_\mathfrak{p}|k_\mathfrak{p})$ are embedable into $G$. Then there exists a Galois extension $K|k$ with Galois group isomorphic to $G$, which for the primes $\mathfrak{p} \in S$ has the given extensions $K_\mathfrak{p}|k_\mathfrak{p}$ as completions.*

The same problem for groups of arbitrary order remains unsolved, we have only the weaker result of Šafarevič which uses an entire different technique.

Let us shortly explain the ideas behind both approaches, the Šafarevič and the Scholz-Reichardt method. By definition, a solvable group is built up by successive extensions of abelian groups. Constructing the required global extension inductively step by step in an abelian way, the first step is given by the theorem of Grunwald-Wang. In the second, and every further step, we have to solve embedding problems with abelian kernel. These are not always solvable, in fact we can reach a dead lock very soon as the following example might indicate (for a proof see [16], Thm. 1.2.4):

*Suppose that $k$ is a field of characteristic $p \neq 2$. Then the quadratic extension $k(\sqrt{a})|k$ can be embedded into a cyclic extension of degree 4 if and only if $a$ is a sum of two squares in $k$.*

The above example shows that, although we did not put local conditions, there might be a global *arithmetic* obstruction to the existence of a solution of our embedding problem. Therefore it is not very promising to solve the embedding problems of every induction step separately. Exactly at this point the methods of Scholz-Reichardt and of Šafarevič differ.

In the Scholz-Reichardt approach (which is not always applicable, see above) one chooses the solutions of the inductively given embedding problems in a very special way, in order to avoid dead locks. Moreover, as Neukirch has shown, one can choose the inductive solutions in such a way, that they satisfy given local conditions at finitely many places.

Šafarevič uses the same special kind of solutions of the inductively given embedding problems ("Scholz solutions"), but in the general situation one can run into a dead lock. The key idea of Šafarevič's approach is to modify the already found solutions of the first $i-1$ induction steps, in order to leave a dead lock within the $i$-th step. This happens in a rather complicated way within a *shrinking procedure*. This method works without restrictive assumptions on the group, but unfortunately it seems to be impossible, to let given local conditions be untouched within the shrinking procedure. Therefore, as the best result of Šafarevič's approach, one can say that every solvable group occurs as a Galois group over $k$ such that the associated local extensions are of a very special type. But we cannot realize given local conditions.

Let us give an outline of our proof of theorem 1 which follows Šafarevič's method:

A first reduction reduces the problem to the solution of split embedding problems with nilpotent kernel. At this point we prove a little bit more than required, namely that a so-called "Scholz-solution" exists, satisfying additional local conditions. (This turns out to be necessary within an induction process.) Let us consider the split embedding problem of finite groups with a $p$-group $N$ as kernel

$$\begin{array}{c} G_k \\ \downarrow \\ 1 \longrightarrow N \longrightarrow E \rightleftarrows G \longrightarrow 1 \,. \end{array}$$

In order to solve these split embedding problems it suffices to consider the generic kernel, i.e. $N = \mathcal{F}(n)/\mathcal{F}(n)^{(\nu)}$ where $\mathcal{F}(n)$ is a free pro-$p$-$G$ operator group of rank $n$ and $\mathcal{F}(n)^{(\nu)}$ denotes the $\nu$-th term of a filtration of $\mathcal{F}(n)$ which we will define below and which refines the descending $p$-central series. [2]

We proceed by induction on $\nu$ whereas $n$ is arbitrary. The case $\nu = 1$ is trivial. The problems which have to be solved within every induction step are of the form

$$\begin{array}{c} G_k \\ \downarrow \\ 1 \to \mathcal{F}(n)^{(\nu)}/\mathcal{F}(n)^{(\nu+1)} \to \mathcal{F}(n)/\mathcal{F}(n)^{(\nu+1)} \rtimes G \to \mathcal{F}(n)/\mathcal{F}(n)^{(\nu)} \rtimes G \to 1. \end{array}$$

This induction step is proved in four substeps: In the first step one shows that this problem is locally solvable everywhere, i.e. for every prime $\mathfrak{p}$ of $k$ the induced local problem

$$\begin{array}{c} G_{k_\mathfrak{p}} \\ \downarrow \\ 1 \longrightarrow \mathcal{F}(n)^{(\nu)}/\mathcal{F}(n)^{(\nu+1)} \longrightarrow E_\mathfrak{p} \longrightarrow (\mathcal{F}(n)/\mathcal{F}(n)^{(\nu)} \rtimes G)_\mathfrak{p} \longrightarrow 1 \end{array}$$

has a solution. This can be done, if the old solution $N_\nu^n|k$ with $G(N_\nu^n|k) = \mathcal{F}(n)/\mathcal{F}(n)^{(\nu)} \rtimes G$ is locally of certain type.

In the second step one uses the local-global principle in order to show that a global solution exists.

In step three and four we modify the global solution, in order to get a proper Scholz-solution, such that the new local problems for the next induction step $\nu + 1$ will be solvable.

---

[2] This refinement, which was proposed by Šafarevič in his correction note, is necessary in order to deal with the case $p = 2$.

Within the induction steps, obstructions to the existence of solutions of the given embedding problems occur. These obstructions really exist and are non-trivial. Šafarevič's idea how to overcome this problem is the following:

We revise the solution found in the $(\nu - 1)$-th induction step. Instead we use the $(\nu - 1)$-th induction step not for $n$ but for a very large $m > n$. The solution of that problem (which exists by the induction hypothesis) induces a new solution for our starting problem via any surjective $G$-invariant homomorphism $\psi : \mathcal{F}(m) \twoheadrightarrow \mathcal{F}(n)$. The shrinking lemma, which we will explain in a moment, tells us that for large enough $m$, we can find $\psi$ in such a way, that all obstructions for the embedding problem in the $\nu$-th induction step for $\mathcal{F}(n)$ vanish. In this way one proves the induction step from $\nu - 1$ to $\nu$ and for every $n$. Having the result for all $n$, we can use the shrinking procedure in the next induction step again.

2. The shrinking procedure

The next proposition is the pure technical skeleton of the shrinking lemma. A similar statement is already contained in the original paper [14].

**Proposition 2** *Let $G$ be a finite group. Suppose that $M$ and $N$ are finitely generated $\mathbb{F}_p[G]$-modules and let $s$, $t \in \mathbb{N}$. Then for large enough[3] $r \in \mathbb{N}$ the following holds:*

*For given elements*
$$z_1, \ldots, z_t \in (\bigoplus_r M)^{\otimes s} \otimes N$$
*there exist $a = (a_i)_{i=1,\ldots,r} \in \mathbb{F}_p^r$, such that*
$$\varphi_a : \bigoplus_r M \longrightarrow M, \quad (x_i)_{i=1,\ldots,r} \longmapsto \sum_{i=1}^r a_i x_i$$
*is a surjective $\mathbb{F}_p[G]$-homomorphism (i.e. not all $a_i$ are zero) and such that the induced homomorphism*
$$\psi_a = (\varphi_a^{\otimes s}) \otimes id : \quad (\bigoplus_r M)^{\otimes s} \otimes N \longrightarrow M^{\otimes s} \otimes N$$
*maps all $z_i$, $i = 1, \ldots, t$, to zero.*

**Proof:**[4] Let $r$ be arbitrary. Then $\varphi_a$ is $\mathbb{F}_p[G]$-invariant and surjective if $a \neq 0$. Let $n = t \dim_{\mathbb{F}_p}(M^{\otimes s} \otimes N)$ and suppose that $r > s\,n$. Consider the set
$$V = \{a \in \mathbb{F}_p^r \,|\, \psi_a(z_1) = \ldots = \psi_a(z_t) = 0\}\,.$$

---
[3]i.e. for all $r \geq r_0$

[4]This proof is based on an idea of *J.Sonn*. We also want to thank *A. Deitmar*.

Then $V$ is the set of common zeroes of $n$ polynomials of degree $s$. It contains the trivial element and by the theorem of Chevalley-Warning (see [15] Chap.1,§2, Thm. 3)) it follows that it also must contain a non-trivial point $a$. Then $\varphi_a$ has all desired properties. □

Let us recall some filtrations:

**Definition 3** *Let $P$ be a pro-$p$-group. Then the descending $p$-central series of $P$ is the filtration $\{P^i\}_{i \geq 1}$ recurrently defined by*

$$P = P^1, \quad P^{i+1} = (P^i)^p[P^i, P],$$

*and the descending central series $\{P_i\}_{i \geq 1}$ of $P$ is given by*

$$P = P_1, \quad P_{i+1} = [P_i, P].$$

*Furthermore, we define a refinement $\{P^{(i,j)}\}$, $i, j \geq 1$, of the descending $p$-central series of $P$ by*

$$P^{(i,j)} := (P^i \cap P_j)P^{i+1}.$$

Obviously we have

$$P^{(i,1)} = P^i \text{ and } P^{(i,j)} = P^{i+1} \text{ for } j > i \geq 1.$$

We introduce the following *notational convention*:

The letter $\nu$ always stands for a pair $(i, j)$, $i \geq j \geq 1$, and we order these pairs lexicographically. We say that $\nu + 1 = (i, j + 1)$ if $i > j$ and $\nu + 1 = (i + 1, 1)$ if $\nu = (i, i)$.

**Lemma 4** (i) *Let $P$ be a pro-$p$-group and let $x \in P^i$, $y \in P^j$ and $a \in p^r \mathbb{Z}_p$. Then*

$(xy)^a \equiv x^a y^a (y, x)^{\binom{a}{2}} \mod P^{i+j+\max\{1,r\}}$,
$(x^a, y) \equiv (x, y)^a ((x, y), x)^{\binom{a}{2}} \mod P^{i+j+1+\max\{1,r\}}$,
$(x, y^a) \equiv (x, y)^a ((x, y), y)^{\binom{a}{2}} \mod P^{i+j+1+\max\{1,r\}}$.

(ii) *Let $P = F$ be a free pro-$p$-group. Then*

$$F^i \cap F_j = (F_j)^{p^{i-j}} \cdot (F_{j+1})^{p^{i-j-1}} \cdot \cdots \cdot F_i$$

*for all $i \geq j \geq 1$.*

**Proof:** All statements are well known and easily proved using induction. For the second statement one uses the fact that $F_j/F_{j+1}$ is a free $\mathbb{Z}_p$-module, which follows from the corresponding theorem of Witt [18] for discrete (finitely generated) free groups. □

**Proposition 5** *For every $\nu = (i,j)$ the $\mathbb{F}_p$-vector space homomorphism*

$$\psi_\nu : \quad (P/P^2)^{\otimes j} \longrightarrow P^{(\nu)}/P^{(\nu+1)}$$
$$\bar{x}_1 \otimes \cdots \otimes \bar{x}_j \longmapsto ([x_1,[x_2,[\cdots,x_j]]\cdots])^{p^{i-j}} \mod P^{(\nu+1)}$$

*is well-defined and surjective.*

**Proof:** Let $S = \{x_\alpha\}_{\alpha \in I}$ be a minimal system of generators of $P$. Then its image $\bar{S} := S \mod P^p[P,P]$ in $P/P^2 = P/P^p[P,P]$ is an $\mathbb{F}_p$-basis of $P/P^2$. Hence the tensors $\bar{x}_{\alpha_1} \otimes \cdots \otimes \bar{x}_{\alpha_j}$, $\alpha_1, \ldots, \alpha_j \in I$ define a basis of $(P/P^2)^{\otimes j}$. Regarding the definition of the $(\nu)$-filtration it follows from lemma 4(ii) (for free groups and then for all pro-$p$-groups), that the elements

$$([x_{\alpha_1},[x_{\alpha_2},[\cdots,x_{\alpha_j}]]\cdots])^{p^{i-j}}, \qquad \alpha_1, \ldots, \alpha_j \in I$$

generate $P^{(\nu)}$ modulo $P^{(\nu+1)}$. Finally an inductive application of lemma 4(i) shows that $\psi_\nu$ is well-defined. □

Recall the definition of a free pro-$p$-$G$ operator group: If $G$ is a finite group and $F_d$ a free pro-$p$-group of rank $d$, then we set

$$\mathcal{F}(d) = \underset{G}{*} F_d$$

for a free pro-$p$-$G$ operator group of rank $d$. The group $\mathcal{F}(d)$ is a free object in the category of pro-$p$-groups with a continuous $G$-action. As pro-$p$-group $\mathcal{F}(d)$ is free of rank $\#G \cdot d$ and one can choose a basis $x_{i,g}$, $i = 1, \ldots, d$, $g \in G$ of $\mathcal{F}(d)$ such that $G$ acts by $g' x_{i,g} = x_{i,g'g}$.

Now we will apply the shrinking process to cohomology groups with respect to $G$ and $\mathcal{F}(d)$ in order to annihilate given cohomology classes.

**Proposition 6** *Let $G$ be a finite group and let $\mathcal{F}(d)$ be the free pro-$p$-$G$ operator group of rank $d$. Let $n, t \in \mathbb{N}$, $k \in \mathbb{Z}$, a finitely generated $\mathbb{F}_p[G]$-module $T$ and $\nu = (i,j)$ be given. Then there exists an $m_0 \geq n$ such that for all $m \geq m_0$ the following holds:*

*For given elements*

$$x_1, \ldots, x_t \in \hat{H}^k(G, \mathcal{F}(m)^{(\nu)}/\mathcal{F}(m)^{(\nu+1)} \otimes T)$$

*there exists a surjective pro-$p$-$G$ operator homomorphism*

$$\psi : \mathcal{F}(m) \twoheadrightarrow \mathcal{F}(n)$$

*such that the induced homomorphism*

$$\psi_* : \hat{H}^k(G, \mathcal{F}(m)^{(\nu)}/\mathcal{F}(m)^{(\nu+1)} \otimes T) \longrightarrow \hat{H}^k(G, \mathcal{F}(n)^{(\nu)}/\mathcal{F}(n)^{(\nu+1)} \otimes T)$$

*maps $x_1, \ldots, x_t$ to zero.*

**Proof:** We set $\mathcal{E}(n,\nu) = \mathcal{F}(n)^{(\nu)}/\mathcal{F}(n)^{(\nu+1)}$. Using dimension shifting we have an isomorphism

$$\hat{H}^k(G, \mathcal{E}(n,\nu) \otimes T) \xrightarrow{\sim} \hat{H}^{-1}(G, \mathcal{E}(n,\nu) \otimes T \otimes A_k),$$

where $A_k = I_G^{\otimes -(k+1)}$. Here $I_G$ denotes the augmentation ideal of $\mathbb{F}_p[G]$ and we have set $I_G^{\otimes 0} = \mathbb{F}_p$ and $I_G^{\otimes -s} = \text{Hom}(I_G^{\otimes s}, \mathbb{F}_p)$ for $s \in \mathbb{N}$.

Since $T$ was arbitrary, we may restrict to the case $k = -1$.

Observe, that for every $d \geq 1$ the canonical surjective $\mathbb{F}_p$-homomorphism defined in lemma 5

$$\theta_\nu(d): \quad (\mathcal{F}(d)/\mathcal{F}(d)^2)^{\otimes j} \twoheadrightarrow \mathcal{F}(d)^{(\nu)}/\mathcal{F}(d)^{(\nu+1)}.$$

is obviously $G$-invariant. For a pro-$p$-$G$ operator homomorphism $\psi: \mathcal{F}(d') \to \mathcal{F}(d)$ we have the compatibility $\theta_\nu(d')\psi_* = \psi_*\theta_\nu(d)$.

Now let $m = rn$, $r$ large enough, such that we can apply proposition 2 with the $G$-module $T$: For given elements

$$z_1, \ldots, z_t \in (\mathcal{F}(m)/\mathcal{F}(m)^2)^{\otimes j} \otimes T = (\mathbb{F}_p[G]^{rn})^{\otimes j} \otimes T,$$

there exists a surjective $\mathbb{F}_p[G]$-homomorphism

$$\bar{\psi}: \mathcal{F}(m)/\mathcal{F}(m)^2 \twoheadrightarrow \mathcal{F}(n)/\mathcal{F}(n)^2,$$

such that $(\bar{\psi}^{\otimes j} \otimes id)(z_\alpha) = 0$, $\alpha = 1, \ldots, t$. By the universal property of free pro-$p$-$G$ operator groups $\bar{\psi}$ extends to a pro-$p$-$G$ operator homomorphism

$$\psi: \mathcal{F}(m) \twoheadrightarrow \mathcal{F}(n),$$

which is necessarily surjective (by the Frattini argument).

Now we consider the commutative diagram

$$\begin{array}{ccccc}
(\mathcal{F}(m)/\mathcal{F}(m)^2)^{\otimes j} \otimes T & \longrightarrow & (\mathcal{E}(m,\nu) \otimes T)_G & \longleftarrow & \hat{H}^{-1}(G, \mathcal{E}(m,\nu) \otimes T) \\
\downarrow \bar{\psi}^{\otimes j} \otimes id & & \downarrow \psi_* \otimes id & & \downarrow \psi_* \\
(\mathcal{F}(n)/\mathcal{F}(n)^2)^{\otimes j} \otimes T & \longrightarrow & (\mathcal{E}(n,\nu) \otimes T)_G & \longleftarrow & \hat{H}^{-1}(G, \mathcal{E}(n,\nu) \otimes T)
\end{array}$$

and choose $z_\alpha$ as a pre-image of the image of $x_\alpha$ in the group $(\mathcal{E}(m,\nu) \otimes T)_G$. Choosing an appropriate $\psi$, the diagram shows that $\psi_*(x_\alpha) = 0$ for $\alpha = 1, \ldots, t$. This proves the proposition. $\square$

We will apply proposition 6 only for $k = 2$ and $k = -2$ and for $k = -2$ we also need the following variant of it, which goes back to an idea of V.V. IšHANOV.

**Proposition 7** Let $G$ be a finite group and let $\mathcal{F}(d)$ be the free pro-$p$-$G$ operator group of rank $d$. Let $n, t \in \mathbb{N}$, a finitely generated $\mathbb{F}_p[G]$-module $T$ and $\nu = (i, j)$ be given. Then for large enough $m \geq n$ the following holds:

(i) For given elements
$$x_1, \ldots, x_t \in H^{-2}(\mathcal{F}(m)/\mathcal{F}(m)^{(\nu)} \rtimes G, \mathcal{F}(m)^{(\nu)}/\mathcal{F}(m)^{(\nu+1)} \otimes T)$$
there exists a surjective pro-$p$-$G$ operator homomorphism
$$\psi : \mathcal{F}(m) \twoheadrightarrow \mathcal{F}(n),$$
such that the induced homomorphism
$$\psi_* : H^{-2}(\mathcal{F}(m)/\mathcal{F}(m)^{(\nu)} \rtimes G, \mathcal{F}(m)^{(\nu)}/\mathcal{F}(m)^{(\nu+1)} \otimes T) \longrightarrow$$
$$H^{-2}(\mathcal{F}(n)/\mathcal{F}(n)^{(\nu)} \rtimes G, \mathcal{F}(n)^{(\nu)}/\mathcal{F}(n)^{(\nu+1)} \otimes T)$$
maps $x_1, \ldots, x_t$ to zero.

(ii) For given elements
$$x_1, \ldots, x_t \in H^{-2}(\mathcal{F}(m)/\mathcal{F}(m)^{(\nu)}, \mathcal{F}(m)^{(\nu)}/\mathcal{F}(m)^{(\nu+1)} \otimes T)$$
there exists a surjective pro-$p$-$G$ operator homomorphism
$$\psi : \mathcal{F}(m) \twoheadrightarrow \mathcal{F}(n),$$
such that the induced homomorphism
$$\psi_* : H^{-2}(\mathcal{F}(m)/\mathcal{F}(m)^{(\nu)}, \mathcal{F}(m)^{(\nu)}/\mathcal{F}(m)^{(\nu+1)} \otimes T) \longrightarrow$$
$$H^{-2}(\mathcal{F}(n)/\mathcal{F}(n)^{(\nu)}, \mathcal{F}(n)^{(\nu)}/\mathcal{F}(n)^{(\nu+1)} \otimes T)$$
maps $x_1, \ldots, x_t$ to zero.

**Proof:** We keep the notations
$$\mathcal{F}(n)/\nu = \mathcal{F}(n)/\mathcal{F}(n)^{(\nu)} \text{ and } \mathcal{E}(n, \nu) = \mathcal{F}(n)^{(\nu)}/\mathcal{F}(n)^{(\nu+1)}.$$

If $\nu = 1$, then the statement to prove is just a special case of proposition 6. So we may assume that $\nu = (i, j) \geq (2, 1)$. Recall that $H^{-2} = H_1$ and consider the exact sequence
$$H_1(\mathcal{F}(n)/\nu, \mathcal{E}(n, \nu) \otimes T) \longrightarrow H_1(\mathcal{F}(n)/\nu \rtimes G, \mathcal{E}(n, \nu) \otimes T) \longrightarrow$$
$$H_1(G, \mathcal{E}(n, \nu) \otimes T) \longrightarrow 0,$$
which is induced by the homological Hochschild-Serre sequence. Since $\mathcal{E}(n, \nu)$ is a trivial $\mathcal{F}(n)/\nu$-module, the universal coefficient formula yields:
$$H_1(\mathcal{F}(n)/\nu, \mathcal{E}(n, \nu) \otimes T) \cong \mathcal{F}(n)/\mathcal{F}(n)^2 \otimes \mathcal{E}(n, \nu) \otimes T$$

and proposition 5 implies the existence of a $G$-invariant surjection
$$(\mathcal{F}(n)/\mathcal{F}(n)^2)^{\otimes(j+1)} \otimes T \twoheadrightarrow H_1(\mathcal{F}(n)/\nu, \mathcal{E}(n, \nu) \otimes T),$$

where $\nu = (i,j)$. This is obviously true for arbitrary $n$, and the maps are compatible. Let $r \geq n$ be the number which has the property, that arbitrary $t$ elements in $(\mathcal{F}(r)/\mathcal{F}(r)^2)^{\otimes(j+1)} \otimes T$ are annihilated by the homomorphism induced by a suitable chosen $G$-invariant surjection $\mathcal{F}(r) \twoheadrightarrow \mathcal{F}(n)$ (and which exists by proposition 2). The above surjection shows that $r$ has the same property with respect to $t$ arbitrary given elements in

$$H_1(\mathcal{F}(r)/\nu, \mathcal{E}(r,\nu) \otimes T).$$

This proves (ii). In order to show (i), let $m \geq r$ be the number, which has the property that arbitrary $t$ elements in

$$H_1(G, \mathcal{E}(m,\nu) \otimes T)$$

are annihilated by the homomorphism induced by a suitable chosen $G$-invariant surjection $\mathcal{F}(m) \twoheadrightarrow \mathcal{F}(r)$ (and which exists by proposition 6). We obtain the commutative and exact diagram

$$\begin{array}{ccccc}
H_1(\mathcal{F}(m)/\nu, \mathcal{E}(m,\nu) \otimes T) & \longrightarrow & H_1(\mathcal{F}(m)/\nu \rtimes G, \mathcal{E}(m,\nu) \otimes T) & \stackrel{\alpha}{\twoheadrightarrow} & H_1(G, \mathcal{E}(m,\nu) \otimes T) \\
\downarrow & & \downarrow \pi & & \downarrow \\
H_1(\mathcal{F}(r)/\nu, \mathcal{E}(r,\nu) \otimes T) & \stackrel{\beta}{\longrightarrow} & H_1(\mathcal{F}(r)/\nu \rtimes G, \mathcal{E}(r,\nu) \otimes T) & \longrightarrow & H_1(G, \mathcal{E}(r,\nu) \otimes T) \\
\downarrow \varepsilon & & \downarrow & & \downarrow \\
H_1(\mathcal{F}(n)/\nu, \mathcal{E}(n,\nu) \otimes T) & \longrightarrow & H_1(\mathcal{F}(n)/\nu \rtimes G, \mathcal{E}(n,\nu) \otimes T) & \longrightarrow & H_1(G, \mathcal{E}(n,\nu) \otimes T)
\end{array}$$

in which the vertical maps are induced by $G$-invariant surjections

$$\mathcal{F}(m) \twoheadrightarrow \mathcal{F}(r) \twoheadrightarrow \mathcal{F}(n),$$

which we choose in the following way:

Let arbitrary elements $x_1, \ldots, x_t \in H_1(\mathcal{F}(m)/\nu \rtimes G, \mathcal{E}(m,\nu) \otimes T)$ be given. Choose $\mathcal{F}(m) \twoheadrightarrow \mathcal{F}(r)$ such that the induced homomorphism annihilates the elements $\alpha(x_1), \ldots, \alpha(x_t)$. Hence $\pi(x_1), \ldots, \pi(x_t)$ are contained in the image of $\beta$ and we choose the surjection $\mathcal{F}(r) \twoheadrightarrow \mathcal{F}(n)$ such that $\varepsilon$ annihilates arbitrary chosen $\beta$-preimages of $\pi(x_1), \ldots, \pi(x_t)$. The composite $\mathcal{F}(m) \twoheadrightarrow \mathcal{F}(n)$ has the desired property. $\square$

3. Some related duality statements

We keep the assumption that $k$ is a global field and introduce some further notations.

**Definition 8** *If $\Omega|k$ is a finite extension of $k$, then we denote by*
$$\delta_\Omega(S) = \delta_\Omega(S(\Omega))$$
*the Dirichlet density of the set $S(\Omega)$ of primes of $\Omega$ given by all extensions of $S = S(k)$ to $\Omega$. For sets $S_1$ and $S_2$ of primes we use the notation*
$$S_1 \subsetneq S_2 :\Longleftrightarrow \delta(S_1\backslash S_2) = 0,$$
*i.e. $S_1$ is contained in $S_2$ up to a set of primes of density zero.*
*Furthermore, we set for a finite Galois extension $\Omega|k$*
$$\begin{aligned} cs(\Omega|k) &:= \{\mathfrak{p} \text{ a prime of } k \mid \mathfrak{p} \text{ splits completely in } \Omega|k\}, \\ \mathrm{Ram}(\Omega|k) &:= \{\mathfrak{p} \text{ a prime of } k \mid \mathfrak{p} \text{ ramifies in } \Omega|k\}. \end{aligned}$$

Let $T \subseteq S$ be sets of primes of $k$ where $S$ is non-empty and contains the archimedean primes $S_\infty$ in the number field case. Let $A$ be a finite $G_S(k)$-module whose order $\#A$ is invertible in $\mathcal{O}_{k,S}$ and let $A' = \mathrm{Hom}(A, \mathcal{O}_S^\times)$.

**Definition 9** *We define the groups $\mathrm{III}^1(k_S, T, A)$ and $\mathrm{coker}(k_S, T, A)$ by the exact sequence*
$$\mathrm{III}^1(k_S, T, A) \hookrightarrow H^1(k_S|k, A) \xrightarrow{res} \prod_T H^1(k_\mathfrak{p}, A) \twoheadrightarrow \mathrm{coker}(k_S, T, A),$$
*where $\prod$ denotes the restricted product with respect to the subgroups $H^1_{nr}(k_\mathfrak{p}, A)$.*

If the Dirichlet density of the set $T$ is equal to $1$ and $A$ is a trivial $G_S(k)$-module, then the *Hasse-principle* holds, i.e.
$$\mathrm{III}^1(k_S, T, A) = 0.$$
This is well known for number fields, [8] VII 13.6, and for functions fields it can be found in [2] Ch.5 §4.

Assume that $T$ is finite, then so is $\mathrm{coker}(k_S, T, A)$. From the local and global duality theorems we obtain the commutative and exact diagram

$$\begin{array}{ccccc}
\mathrm{III}^1(k_S, S\backslash T, A') & \hookrightarrow & H^1(k_S|k, A') & \to & \prod_{S\backslash T} H^1(k_\mathfrak{p}, A') \\
\uparrow & & \parallel & & \uparrow \\
\mathrm{III}^1(k_S, A') & \hookrightarrow & H^1(k_S|k, A') & \to & \prod_S H^1(k_\mathfrak{p}, A') & \longrightarrow & H^1(k_S|k, A)^* \\
& & & & \uparrow & & \uparrow \\
& & & & \prod_T H^1(k_\mathfrak{p}, A') & \xrightarrow{\sim} & \prod_T H^1(k_\mathfrak{p}, A)^* \\
& & & & & & \uparrow \\
& & & & & & \mathrm{coker}(k_S, T, A)^*.
\end{array}$$

This diagram implies the

**Lemma 10** *Assume that $T$ is finite, then there is a canonical exact sequence*

$$0 \to \text{III}^1(k_S, A') \to \text{III}^1(k_S, S\backslash T, A') \to \text{coker}(k_S, T, A)^* \to 0.$$

Under a mild restriction, the Hasse-principle also holds for the module $\mu_m$. Via lemma 10 one can deduce the theorem of *Grunwald-Wang*:

$$\text{coker}(k_S, T, A) = 0,$$

if the density $\delta(S)$ is equal to 1, $A$ is a trivial $G_S(k)$-module and $T$ is finite not containing primes above 2 in the number field case (the last assumption guarantees that we avoid the so-called *special case*).

We have the following application of the Grunwald-Wang theorem to embedding problems with induced $G$-modules as kernel.

**Proposition 11** *Let $K|k$ be a finite Galois extension of global fields with Galois group $G = G(K|k)$ and let $A = \mathbb{F}_p[G]^n$. Then the embedding problem*

$$\begin{array}{ccccccccc}
& & & & & & G_k & & \\
& & & & & & \downarrow & & \\
1 & \to & A & \to & E & \to & G & \to & 1
\end{array}$$

*is properly solvable.*

**Proof:** Since $H^2(G, A) = 0$, the embedding problem has a solution $\psi_0 : G_k \to E$. Let $\mathfrak{p}_1, \ldots, \mathfrak{p}_r$ be primes of $k$ which split completely in $K$ and let $\varphi_i : G_{k_{\mathfrak{p}_i}} \to A$ be homomorphisms such that their images generate $A$. Furthermore we assume that the primes $\mathfrak{p}_i$ do not divide 2 if $k$ is a number field. Then by the Grunwald-Wang theorem the map $res = (res_i)_i$

$$\begin{array}{ccc}
H^1(k, A) & \xrightarrow{res} & \prod_{i=1}^r H^1(k_{\mathfrak{p}_i}, A) = \prod_{i=1}^r \text{Hom}(G_{k_{\mathfrak{p}_i}}, A) \\
\downarrow \wr & & \downarrow \wr \\
H^1(K, \mathbb{F}_p{}^n) & \to & \prod_{i=1}^r \prod_{\mathfrak{P}|\mathfrak{p}_i} H^1(K_{\mathfrak{P}}, \mathbb{F}_p{}^n)
\end{array}$$

is surjective. Let $[x] \in H^1(k, A)$ such that $res_i[x] = \varphi_i - \psi_0|_{G_{k_{\mathfrak{p}_i}}}$ for $i = 1, \ldots, r$. Then $\psi = x \cdot \psi_0 : G_k \to E$ is a new solution of the embedding problem which is proper, since $\psi|G_{k_{\mathfrak{p}_i}} = \varphi_i$ for $i = 1, \ldots, r$, hence $\psi(G_K) = A$. □

We introduce some additional notations. For a prime number $p \neq \text{char}(k)$ we denote the set of primes of $k$ with residue characteristic $p$ by $\boldsymbol{S_p} = \boldsymbol{S_p(k)}$. The set $S_p(k)$ is finite and it is empty if $k$ is a function field. In the number field case we denote by $S_\infty = S_\infty(k)$ the set of archimedean places of $k$. In the function field case we *choose* any finite, non-empty set of primes of $k$ and call it $\boldsymbol{S_\infty = S_\infty(k)}$. For every extension field $K|k$ we denote by $S_\infty(K)$ the set of primes of $K$ which lie over $S_\infty(k)$. For a $G_k$-module $A$ we denote by $k(A)$ the minimal trivializing extension of $k$, i.e. $G_{k(A)}$ is the kernel of the homomorphism $G_k \to \mathrm{Aut}(A)$ given by the action of $G_k$ on $A$. The next technical lemma will be needed later.

**Lemma 12** *Let $k$ be a global field, $p \ne \mathrm{char}(k)$ a prime number and assume that we are given sets of primes of $k$*

$$S' \supseteq S \supseteq T \supseteq S_p \cup S_\infty \,,$$

*where $T$ is finite. Let $A$ be a finite $G_S$-module which is annihilated by $p$. In addition suppose that we are given a finite subextension $N \subseteq k_S$, with*
  a) $k(A) \subseteq N$,
  b) $S' \smallsetminus T \subseteq cs(N|k)$,
  c) $\mu_p \subseteq N$.

*Consider the diagram with solid arrows (the lines are not exact)*

$$\begin{array}{ccccc}
\mathrm{coker}(k_S, T, A) & \hookrightarrow & \mathrm{III}^1(k_S, S \smallsetminus T, A')^* & \xrightarrow{\eta} & H^1(N|k, A')^* \\
\downarrow & & \downarrow \phi & & \| \\
\mathrm{coker}(k_{S'}, T, A) & \hookrightarrow & \mathrm{III}^1(k_{S'}, S' \smallsetminus T, A')^* & \xrightarrow{\eta'} & H^1(N|k, A')^*
\end{array}$$

*in which $A' := A^*(1) = \mathrm{Hom}(A, \mu_p)$. The horizontal maps on the left are induced by lemma 10 and those on the right come from the Hochschild-Serre sequence and from conditions a),b),c).*

*Then in the above situation a natural dotted arrow $\phi$ making the diagram commutative exists. If in addition*

$$cs(N|k) \subsetneq S' \smallsetminus T,$$

*then the surjection $\eta'$ is an isomorphism.*

**Proof:** First observe that the homomorphism $\eta'$ is obtained from the following commutative and exact diagram

$$\begin{array}{ccc}
H^1(N_{S'}, A') & \xrightarrow{\iota} & \prod_{S' \backslash T} H^1(N_{\mathfrak{P}}, A') \\
\uparrow & & \uparrow \\
\mathrm{III}^1(k_{S'}, S'\backslash T, A') \hookrightarrow H^1(k_{S'}, A') & \longrightarrow & \prod_{S' \backslash T} H^1(k_{\mathfrak{p}}, A') \\
\quad {}_{(\eta')^*}\searrow \quad \uparrow \quad \nearrow {}_{0} & & \\
H^1(N|k, A') & &
\end{array}$$

and in a similar way we get the homomorphism $\eta$.

Now consider the exact and commutative diagram with natural homomorphisms:

$$\begin{array}{ccccc}
\mathrm{III}^1(k_S, S\backslash T, A') & \hookrightarrow & H^1(G_S, A') & \longrightarrow & \prod_{S\backslash T} H^1(k_\mathfrak{p}, A') \\
\uparrow \kappa & & \| & & \uparrow \\
\ker(\alpha) & \hookrightarrow & H^1(G_S, A') & \xrightarrow{\alpha} & \prod_{S\backslash T} H^1(k_\mathfrak{p}, A') \times \prod_{S'\backslash S} H^1_{nr}(k_\mathfrak{p}, A') \\
\downarrow \varepsilon & & \downarrow & & \uparrow \\
\mathrm{III}^1(k_{S'}, S'\backslash T, A') & \hookrightarrow & H^1(G_{S'}, A') & \longrightarrow & \prod_{S\backslash T} H^1(k_\mathfrak{p}, A') \times \prod_{S'\backslash S} H^1(k_\mathfrak{p}, A') \\
& & \downarrow & & \downarrow \\
& & H^1(k_{S'}|k_S, A')^{G_S} & \xrightarrow{\beta} & \prod_{S'\backslash S} H^1(T_\mathfrak{p}, A')^{G_{k_\mathfrak{p}}} \,,
\end{array}$$

in which $T_\mathfrak{p} \subseteq G_{k_\mathfrak{p}}$ denotes the inertia group. Observe that $\beta$ is injective, since $A'$ is a trivial $G(k_{S'}|k_S)$-module and since $k_S$ has by definition no extensions in $k_{S'}$ which are unramified at all places in $S' \smallsetminus S$. Diagram chasing then shows that $\varepsilon$ is an isomorphism. We can now define $\phi$ as the dual homomorphism to $\kappa \circ (\varepsilon^{-1})$. Conditions a),b),c) imply, that $H^1(N|k, A')$ is canonically contained in the groups $\mathrm{III}^1(k_S, S\backslash T, A')$, $\ker(\alpha)$ and $\mathrm{III}^1(k_{S'}, S'\backslash T, A')$. Thus we have constructed $\phi$ and we see that the right part of the diagram is commutative. But that the left part of the diagram is also commutative can be seen from the diagram before lemma 10 which defines the occurring maps.

Now assume that $cs(N|k) \subsetneq S' \smallsetminus T$. Then $\delta_N(S'\backslash T) = 1$ and by the Hasse-principle the homomorphism $\iota$ in the exact commutative diagram at the beginning of the proof is injective. Hence the inclusion $(\eta')^*$ is an isomorphism. □

4. Construction of certain cohomology classes

Let $k_\mathfrak{p}$ be a local field and let $A$ be a $G_{k_\mathfrak{p}}$-module. We call a class $x_\mathfrak{p} \in H^1(k_\mathfrak{p}, A)$ **cyclic**, if it is split by a cyclic extension of $k_\mathfrak{p}$, i.e. if there exists a cyclic extension $K_\mathfrak{p}|k_\mathfrak{p}$ such that $x_\mathfrak{p}$ lies in the kernel of the restriction map $H^1(k_\mathfrak{p}, A) \to H^1(K_\mathfrak{p}, A)$. If $A$ is unramified (i.e. the inertia group acts trivially), then we call $x_\mathfrak{p}$ **unramified** if it is contained in the unramified part $H^1_{nr}(k_\mathfrak{p}, A)$ of $H^1(k_\mathfrak{p}, A)$. In particular, if $x_\mathfrak{p}$ is unramified, then it is cyclic.

The following existence theorem is based on Čebotarev's density theorem and will be used in step 4 of the proof of theorem 15 below. In the case $A = \mu_p$, via Kummer theory, theorem 13 is equivalent to Šafarevič's theorem about the existence of certain algebraic numbers ([12]).

**Theorem 13** *Let $p$ be a prime number and let $\Omega|K|k$ be finite Galois extensions of global fields of characteristic different to $p$, where $K$ contains the group $\mu_p$ of $p$-th roots of unity. Let $T$ be a finite set of primes of $k$ containing $\mathrm{Ram}(\Omega|k) \cup S_p \cup S_\infty$[5] and let $S = cs(\Omega|k) \cup T$.*

*Let $A$ be a finite $\mathbb{F}_p[G(K|k)]$-module and assume we are given a class $y$ in $H^1(k_S|K, A)$ such that*

$y_\mathfrak{P}$ *is unramified for* $\mathfrak{P} \in T(K)$ *and* $y_\mathfrak{P} = 0$ *for* $\mathfrak{P} \cap k \in \mathrm{Ram}(K|k) \cup S_p \cup S_\infty$.

*Then there exists an element $x \in H^1(k_S|k, A)$ such that*

$x_\mathfrak{p} = (cor\,_k^K y)_\mathfrak{p}$ *for* $\mathfrak{p} \in T$ *and* $x_\mathfrak{p}$ *is cyclic for all* $\mathfrak{p} \notin T$.

**Proof:** Putting $x = cor\,_k^K z$, it suffices to construct a $z \in H^1(k_S|K, A)$ with
(a) $z_\mathfrak{P} = y_\mathfrak{P}$ for $\mathfrak{P} \in T(K)$,
(b) If $\mathfrak{P} \notin T(K)$ and $z_\mathfrak{P}$ is ramified (i.e. not contained in $H^1_{nr}(K_\mathfrak{P}, A)$), then $z_\mathfrak{P}$ is cyclic and $z_{\sigma\mathfrak{P}} = 0$ for every $\sigma \in G(K|k) \smallsetminus \{1\}$.

We first prove the existence of $z$ in the case $A = \mu_p$, where the cyclicity condition is trivially satisfied. Assume first that $p$ is odd. We will apply the method of [12] in order to obtain the element $z$ which we are looking for as a product of two members of a sequence

$$z_1, z_2, z_3, \ldots \in H^1(k_S|K, \mu_p),$$

which will be constructed having the following properties.
(1) There exists a prime $\mathfrak{P}_i \in S \smallsetminus T(K)$, such that
  · $(z_i)_{\mathfrak{P}_i} \equiv \mathrm{Frob}\,_{\mathfrak{P}_i}$ modulo $H^1_{nr}(K_{\mathfrak{P}_i}, \mu_p)$,
  · $(z_i)_\mathfrak{P} \in H^1_{nr}(K_\mathfrak{P}, \mu_p)$ for all $\mathfrak{P} \neq \mathfrak{P}_i$,
(2) $(z_i)_\mathfrak{P} = \frac{1}{2} y_\mathfrak{P}$ for $\mathfrak{P} \in T(K)$,
(3) $(z_{n+1})_{\sigma\mathfrak{P}_i} = -(z_i)_{\sigma\mathfrak{P}_i}$ for $i \leq n$ and all $\sigma \in G(K|k) \smallsetminus \{1\}$.

In (1) we view $\mathrm{Frob}\,_\mathfrak{P}$ as an element of $H^1(K_\mathfrak{P}, \mu_p)/H^1_{nr}(K_\mathfrak{P}, \mu_p)$ via

$$G(K_\mathfrak{P}^{nr}|K_\mathfrak{P}) \twoheadrightarrow H^1_{nr}(K_\mathfrak{P}, \mathbb{Z}/p\mathbb{Z})^\vee \cong H^1(K_\mathfrak{P}, \mu_p)/H^1_{nr}(K_\mathfrak{P}, \mu_p).$$

where the isomorphism is induced by the local duality theorem, see [17] chap.II §5.2, §5.5. Assume that we have already constructed $z_1, \ldots, z_n$ ($n \geq 0$) and set $T_n = T \cup \{(\mathfrak{P}_1 \cap k), \ldots, (\mathfrak{P}_n \cap k)\}$, i.e. $T_n(K)$ consists of $T(K)$ and of all $G(K|k)$-conjugates of $\mathfrak{P}_1, \ldots, \mathfrak{P}_n$. Observe that $T_n \subseteq S$ and consider the commutative and exact diagram

---
[5]Regard our notational convention concerning $S_p \cup S_\infty$ in the function field case.

$$
\begin{array}{ccc}
H^1(k_S|K,\mu_p) \longrightarrow \prod_{T_n} H^1(K_{\mathfrak{P}},\mu_p) \times \bigoplus_{S\setminus T_n} H^1(K_{\mathfrak{P}},\mu_p)/H^1_{nr} & \longrightarrow & H^1(k_{T_n}|K,\mathbb{Z}/p\mathbb{Z})^{\vee} \\
& \downarrow & \downarrow \\
& \prod_T H^1(K_{\mathfrak{P}},\mu_p) \longrightarrow & H^1(\Omega|K,\mathbb{Z}/p\mathbb{Z})^{\vee} \\
& \uparrow & \uparrow \\
H^1(k_S|K,\mu_p) \longrightarrow \prod_S H^1(K_{\mathfrak{P}},\mu_p) & \longrightarrow & H^1(k_S|K,\mathbb{Z}/p\mathbb{Z})^{\vee},
\end{array}
$$

in which $\prod$ denotes the product of local cohomology groups restricted with respect their subgroups of unramified elements and $^{\vee}$ is the Pontrjagin dual. The lower line is part of the long exact sequence of Tate-Poitou.

Consider the element

$$\xi \in \prod_{T_n} H^1(K_{\mathfrak{P}},\mu_p) \times \bigoplus_{S\setminus T_n} H^1(K_{\mathfrak{P}},\mu_p)/H^1_{nr}$$

given by $\xi_{\mathfrak{P}} = \frac{1}{2}y_{\mathfrak{P}}$ for $\mathfrak{P} \in T(K)$, $\xi_{\sigma\mathfrak{P}_i} = -(z_i)_{\sigma\mathfrak{P}_i}$ for $i \leq n$, $\sigma \in G(K|k)\setminus\{1\}$ and $\xi_{\mathfrak{P}} = 0$ for all other $\mathfrak{P}$. Then $\xi$ has the same image in $H^1(\Omega|K,\mathbb{Z}/p\mathbb{Z})^{\vee}$ as $y$ ($\in H^1(k_S|K,\mu_p)$). By the exactness of the lower line, this image is trivial. Using Čebotarev's density theorem, we can choose a prime $\mathfrak{P}_{n+1} \in S\setminus T_n(K)$ such that the image of $-\xi$ in $H^1(k_{T_n}|K,\mathbb{Z}/p\mathbb{Z})^{\vee}$ is equal to $\mathrm{Frob}\,\mathfrak{P}_{n+1}$. By the exactness of the upper line we find a class $z_{n+1} \in H^1(k_S|K,\mu_p)$ with properties (1),(2),(3).

Let $G(K|k)\setminus\{1\} = \{\sigma_1,\ldots,\sigma_r\}$ and consider for $n = 1,2,\ldots$ the maps

$$\psi_n : \{z_1,\ldots,z_n\} \longrightarrow \prod_1^r \mu_p\,,$$

given by

$$\psi_n(z_i) = \left((z_i)_{\sigma_1\mathfrak{P}_i}(\mathrm{Frob}\,_{\sigma_1\mathfrak{P}_i}),\ldots,(z_i)_{\sigma_r\mathfrak{P}_i}(\mathrm{Frob}\,_{\sigma_r\mathfrak{P}_i})\right).$$

(Observe that by construction $z_i \in H^1_{nr}(K_{\sigma_j\mathfrak{P}_i},\mu_p)$ for $j = 1,\ldots,r$.)
By the shoe box principle, there exists an $N$ with $\psi_N(z_N) = \psi_N(z_i)$ for some $i < N$. We claim that

$$z = z_i + z_N \in H^1(k_S|K,\mu_p)$$

satisfies conditions (a) and (b) above. Indeed, (a) is trivial by condition (2). It therefore remains to show that, if $z_{\mathfrak{P}}$ is ramified for some $\mathfrak{P}$, then $\mathfrak{P} \in S\setminus T(K)$ and $z_{\sigma\mathfrak{P}} = 0$ for $\sigma \in G(K|k)\setminus\{1\}$. By construction, $z$ is only ramified at $\mathfrak{P}_i$ and $\mathfrak{P}_N$ and by condition (3) we know that $z_{\sigma\mathfrak{P}_i} = 0$ for $\sigma \neq 1$. In order to show the corresponding statement at $\sigma\mathfrak{P}_N$ ($\sigma \neq 1$), recall that for arbitrary classes $a,b \in H^1(k_S|K,\mu_p)$ we have the product formula

$$\prod_{\mathfrak{P}\in S(K)}(a,b)_{\mathfrak{P}} = 1,$$

for the Hilbert symbol $((a,b)_{\mathfrak{P}}$ is defined as the image of $a \cup b$ under the trace homomorphism $H^2(K_{\mathfrak{P}}, \mu_p^{\otimes 2}) \xrightarrow{\sim} \mu_p)$. Since $z_i$ and $z_N$ are unramified at $\sigma \mathfrak{P}_N$, it suffices to show that their values on $\operatorname{Frob}_{\sigma\mathfrak{P}_N}$ are mutually inverse in $\mu_p$. We have

$$
\begin{aligned}
z_N(\operatorname{Frob}_{\sigma\mathfrak{P}_N}) &= z_i(\operatorname{Frob}_{\sigma\mathfrak{P}_i}) && \text{because } \psi_n(z_i) = \psi_n(z_N), \\
&= (z_i, \sigma z_i)_{\sigma\mathfrak{P}_i} && \text{by condition (1) for } z_i, \\
&= (z_i, \sigma z_i)_{\mathfrak{P}_i}^{-1} && \text{by the product formula and by (1), (2),} \\
&= (z_i, \sigma z_N)_{\mathfrak{P}_i} && \text{by condition (3),} \\
&= (z_i, \sigma z_N)_{\sigma\mathfrak{P}_N}^{-1} && \text{by the product formula and by (1), (2),} \\
&= z_i(\operatorname{Frob}_{\sigma\mathfrak{P}_N})^{-1} && \text{by condition (1) for } z_N.
\end{aligned}
$$

This finishes the case $A = \mu_p$, $p$ odd. In the case $p = 2$ we have to modify the method and we will obtain $z$ as a product of three other elements. We use the combinatorial method of [5] chap.5 §3.

Let $\{G_1, G_2, G_3\}$ be a partition of the set $G(K|k) \smallsetminus \{1\}$ such that $G_1$ consists of all elements of order 2 and $G_2 = G_3^{-1}$.

We construct recurrently a sequence $z_1, z_2, \ldots$ of elements in $H^1(k_S|K, \mu_2)$ satisfying the following properties

(1) There exists a prime $\mathfrak{P}_i \in S \smallsetminus T(K)$, such that
- $(z_i)_{\mathfrak{P}_i} \equiv \operatorname{Frob}_{\mathfrak{P}_i}$ modulo $H^1_{nr}(K_{\mathfrak{P}_i}, \mu_2)$,
- $(z_i)_{\mathfrak{P}} \in H^1_{nr}(K_{\mathfrak{P}}, \mu_2)$ for all $\mathfrak{P} \neq \mathfrak{P}_i$,

(2) $(z_i)_{\mathfrak{P}} = y_{\mathfrak{P}}$ for $\mathfrak{P} \in T(K)$,

(3) $(z_{n+1})_{\sigma\mathfrak{P}_i} = 0$ for $i \leq n$ and all $\sigma \in G_1$.

(4) If $\psi_n(z_i) \neq \psi_n(z_j)$ for all $j$ with $n \geq j > i$, then
$$(z_{n+1})_{\sigma\mathfrak{P}_i} = \begin{cases} 0 & \text{if } \sigma \in G_2, \\ (\sigma^2 z_i)_{\sigma\mathfrak{P}_i} & \text{if } \sigma \in G_3, \end{cases}$$
and otherwise: $(z_{n+1})_{\sigma\mathfrak{P}_i} = \begin{cases} (z_i)_{\sigma\mathfrak{P}_i} & \text{if } \sigma \in G_2, \\ 0 & \text{if } \sigma \in G_3. \end{cases}$

The existence of the sequence of classes $z_1, z_2, \ldots$ is proved similar to the case of odd $p$. In addition, we obtain

*Claim:* $\quad (z_i)_{\sigma\mathfrak{P}_i} = 0$ for $\sigma \in G_1$.

*Proof of the claim:* Let $\tilde{z}_i \in K^\times$ be a representative of $z_i \in H^1(k_S|K, \mu_2) \subseteq H^1(K, \mu_2) \cong K^\times / K^{\times 2}$. By condition (1) we have $(\tilde{z}_i) = \mathfrak{P}_i \mathfrak{A}^2$ with a fractional ideal $\mathfrak{A}$ of $\mathcal{O}_{K, S_\infty}$. By Čebotarev's density theorem there exists a prime ideal $\mathfrak{Q} \notin T(K)$ of $\mathcal{O}_{K, S_\infty}$ with $\mathfrak{Q} \neq \sigma\mathfrak{Q}$ such that $\mathfrak{Q} = \mathfrak{A} \cdot (x)$ with $x \in K^\times$. Hence in $\mathcal{O}_{K, S_\infty}$ we have $(\tilde{z}_i x^2) = \mathfrak{P}_i \mathfrak{Q}^2$. Thus we may assume (using condition (2)) that $\tilde{z}_i \in \mathcal{O}_{K, S_\infty}$, $v_{\mathfrak{P}_i}(\tilde{z}_i) = 1$, $v_{\mathfrak{P}}(\tilde{z}_i) = 0$ for $\mathfrak{P} \in T(K)$, $(\tilde{z}_i)_{\mathfrak{P}}$ is a square in $K_{\mathfrak{P}}$ for $\mathfrak{P} \in S_\infty$ and $\tilde{z}_i$ and $\sigma\tilde{z}_i$ are coprime in $\mathcal{O}_{K, S_\infty}$.

In addition, choose $\delta \in \mathcal{O}_{K, S_\infty}$ such that $K = K^\sigma(\delta)$ and $\delta^2 \in \mathcal{O}_{K^\sigma, S_\infty}$, where $K^\sigma$ is the fixed field of $K$ with respect to $\langle \sigma \rangle$. Then there are $a, b \in K^\sigma$ with $\tilde{z}_i = a + b\delta$. In particular, $2a \in \mathcal{O}_{K^\sigma, S_\infty}$ and $2b\delta \in \mathcal{O}_{K, S_\infty}$. Let $\Sigma$ be the set of prime divisors of $2b\delta$ which are not in $S_2 \cup S_\infty \cup \mathrm{Ram}(K|K^\sigma)$. We obtain

$$\begin{aligned}
z_i(\mathrm{Frob}_{\sigma\mathfrak{P}_i}) &= (\tilde{z}_i, \sigma\tilde{z}_i)_{\sigma\mathfrak{P}_i} && \text{by definition of } \tilde{z}_i, \\
&= (2b\delta, \sigma\tilde{z}_i)_{\sigma\mathfrak{P}_i} && \text{since } (\tilde{z}_i - \sigma\tilde{z}_i)^{-1}\tilde{z}_i \in U^1_{\sigma\mathfrak{P}_i}, \\
&= \prod_{\mathfrak{P} \neq \sigma\mathfrak{P}_i}(2b\delta, \sigma\tilde{z}_i)_\mathfrak{P} && \text{by the product formula,} \\
&= \prod_{\mathfrak{P} | 2b\delta\infty}(2b\delta, \sigma\tilde{z}_i)_\mathfrak{P} && \text{if } \mathfrak{P}|\sigma\tilde{z}_i \text{ and } \mathfrak{P} \neq \sigma\mathfrak{P}_i, \text{ then} \\
&&& 2|v_\mathfrak{P}(\sigma\tilde{z}_i) \text{ and } v_\mathfrak{P}(2b\delta) = 0, \\
&= \prod_{\mathfrak{P} \in \Sigma}(2b\delta, \sigma\tilde{z}_i)_\mathfrak{P} && \sigma\tilde{z}_i \in K_\mathfrak{P}^{\times 2} \text{ for } \mathfrak{P} \in S_2 \cup S_\infty \cup \mathrm{Ram}(K|k) \\
&&& \text{by condition (2),} \\
&= \prod_{\mathfrak{P} \in \Sigma}(2b\delta, a)_\mathfrak{P} && \sigma\tilde{z}_i = a - b\delta \text{ and } \mathfrak{P} \nmid \sigma\tilde{z}_i \text{ for } \mathfrak{P} \in \Sigma.
\end{aligned}$$

The last product is easily seen to be unity: If $\mathfrak{P} \in \Sigma$ is inert in $K|K^\sigma$, then $a \in K_\mathfrak{P}^{\times 2}$, hence $(2b\delta, a)_\mathfrak{P} = 1$ and if $\mathfrak{P} \in \Sigma$ splits in $K|K^\sigma$, then $\sigma\mathfrak{P} \in \Sigma$ and

$$(2b\delta, a)_\mathfrak{P} \cdot (2b\delta, a)_{\sigma\mathfrak{P}} = (-1, a)_\mathfrak{P} = 1.$$

This proves the claim.

Now choose $N$ minimal, such that

$$\psi_N(z_i) = \psi_N(z_j) = \psi_N(z_N)$$

for numbers $i < j < N$. We claim that $z = z_i + z_j + z_N$ satisfies conditions (a) and (b). Indeed, (a) follows immediately from (2). It therefore remains to show that, if $z_\mathfrak{P}$ is ramified for some $\mathfrak{P}$, then $\mathfrak{P} \in S \smallsetminus T(K)$ and $z_{\sigma\mathfrak{P}} = 0$ for $\sigma \in G(K|k) \smallsetminus \{1\}$. By construction, $z$ is only ramified at $\mathfrak{P}_i, \mathfrak{P}_j$ and $\mathfrak{P}_N$.

For $\sigma \in G_1$ and $N \geq s, t \geq 1$ we have

$$z_s(\mathrm{Frob}_{\sigma\mathfrak{P}_t}) = 1.$$

which is seen by condition (3) for $s > t$, by the claim for $s = t$ and follows for $s < t$ by (1)-(3) and the product formula

$$z_s(\mathrm{Frob}_{\sigma\mathfrak{P}_t}) = (z_s, \sigma z_t)_{\sigma\mathfrak{P}_t} = (\sigma z_s, z_t)_{\mathfrak{P}_t} = (\sigma z_s, z_t)_{\sigma\mathfrak{P}_s} = z_t(\mathrm{Frob}_{\sigma\mathfrak{P}_s}) = 1.$$

Summing up, we obtain $z_{\sigma\mathfrak{P}_i} = z_{\sigma\mathfrak{P}_j} = z_{\sigma\mathfrak{P}_N} = 0$ for $\sigma \in G_1$.

If $\sigma \in G_2$, then by condition (4)

$$z_{\sigma\mathfrak{P}_i} = (z_i)_{\sigma\mathfrak{P}_i} + (z_j)_{\sigma\mathfrak{P}_i} + (z_N)_{\sigma\mathfrak{P}_i} = (z_i)_{\sigma\mathfrak{P}_i} + 0 + (z_i)_{\sigma\mathfrak{P}_i} = 0.$$

Furthermore, since $(z_N)_{\sigma\mathfrak{P}_j} = 0$ by condition (4) and

$$\begin{aligned}
z_i(\mathrm{Frob}_{\sigma\mathfrak{P}_j}) &= (z_i, \sigma z_j)_{\sigma\mathfrak{P}_j} && \text{by condition (1),} \\
&= (z_i, \sigma z_j)_{\mathfrak{P}_i} && \text{by the product formula,} \\
&= (\sigma^{-1} z_i, z_j)_{\sigma^{-1}\mathfrak{P}_i} && \text{using Galois invariance,} \\
&= z_j(\mathrm{Frob}_{\sigma^{-1}\mathfrak{P}_i}) && \text{by condition (1),} \\
&= \sigma^{-2} z_i(\mathrm{Frob}_{\sigma^{-1}\mathfrak{P}_i}) && \text{by condition (4),} \\
&= z_i(\mathrm{Frob}_{\sigma\mathfrak{P}_i}) && \text{using Galois invariance,} \\
&= z_j(\mathrm{Frob}_{\sigma\mathfrak{P}_j}) && \text{because } \psi_N(z_i) = \psi_N(z_j)\,,
\end{aligned}$$

we obtain $z_{\sigma\mathfrak{P}_j} = 0$ for $\sigma \in G_2$. Similar calculations show that $z_{\sigma\mathfrak{P}_N} = 0$ for $\sigma \in G_2$ and in the same way one verifies that the local classes $z_{\sigma\mathfrak{P}_i}$, $z_{\sigma\mathfrak{P}_j}$ and $z_{\sigma\mathfrak{P}_N}$ also vanish if $\sigma \in G_3$. This finishes the proof for $A = \mu_p$.

The general case will be proven by induction on $\dim_{\mathbb{F}_p} A$. Let $A = A' \oplus \mu_p$. For each $z \in H^1(K_S|K, A)$ let $z = z' + z''$ be the decomposition of $z$ into the components $z' \in H^1(K_S|K, A')$ and $z'' \in H^1(K_S|K, \mu_p)$, and similarly $z_\mathfrak{P} = z'_\mathfrak{P} + z''_\mathfrak{P}$ for $z_\mathfrak{P} \in H^1(K_\mathfrak{P}, A)$. By induction we find an element $z' \in H^1(K_S|K, A')$ such that

(a') $z'_\mathfrak{P} = y'_\mathfrak{P}$ for $\mathfrak{P} \in T(K)$.
(b') If $\mathfrak{P} \notin T(K)$ and $z'_\mathfrak{P}$ is ramified, then $z'_\mathfrak{P}$ is cyclic and $z'_{\sigma\mathfrak{P}} = 0$ for every $\sigma \in G(K|k) \smallsetminus \{1\}$.

Let $K'|K$ be the extension defined by the homomorphism $z'$ (i.e. $z' : G(k_S|K) \twoheadrightarrow G(K'|K) \subseteq A'$) and let $\tilde{K}$ be its Galois closure over $k$. By construction, $K'|K$ and hence also $\tilde{K}|K$ is unramified at all primes in $T(K)$ and it only ramifies at primes in $cs(\Omega|k)(K)$.

Set $T' = T \cup \mathrm{Ram}(\tilde{K}|k)$, $\Omega' = \tilde{K}\Omega$, $S' = cs(\Omega'|k) \cup T'$ and suppose we have found a class $\tilde{y} \in H^1(K_{S'}|K, \mu_p)$ with
- $\tilde{y}_\mathfrak{P} = y''_\mathfrak{P}$ for $\mathfrak{P} \in T(K)$,
- $\tilde{y}_\mathfrak{P} = 0$ for $\mathfrak{P} \in T' \smallsetminus T(K)$.

Then we can apply the induction hypothesis to the field extensions $\Omega'|K|k$, the sets $T'$, $S'$ and the module $A = \mu_p$ in order to find a class $z'' \in H^1(K_{S'}|K, \mu_p)$ with

(a'') $z''_\mathfrak{P} = \tilde{y}_\mathfrak{P}$ for $\mathfrak{P} \in T'(K)$.
(b'') If $\mathfrak{P} \notin T'(K)$ and $z''_\mathfrak{P}$ is ramified, then $z''_\mathfrak{P}$ is cyclic and $z''_{\sigma\mathfrak{P}} = 0$ for every $\sigma \in G(K|k) \smallsetminus \{1\}$.

Regarding $T' \smallsetminus T \subseteq cs(\Omega|k)$, hence $S' \subseteq S$, it is now easily verified, that the class $z = z' + \inf z'' \in H^1(K_S|K, A)$ satisfies conditions (a) and (b). Indeed, for $\mathfrak{P} \in T(K)$ we have
$$z_\mathfrak{P} = z'_\mathfrak{P} + z''_\mathfrak{P} = y'_\mathfrak{P} + y''_\mathfrak{P} = y_\mathfrak{P}\,.$$

For $\mathfrak{P} \in T'\backslash T$ we get
$$z_\mathfrak{P} = z'_\mathfrak{P} + z''_\mathfrak{P} = z'_\mathfrak{P}\,.$$

Therefore $z_\mathfrak{P}$ is cyclic and if $z_\mathfrak{P} = z'_\mathfrak{P}$ is ramified, then the underlying prime $\mathfrak{p}$ of $\mathfrak{P}$ splits completely in $\Omega$ and $z_{\sigma\mathfrak{P}} = z'_{\sigma\mathfrak{P}} = 0$ for $\sigma \in G(K|k)\setminus\{1\}$. Let finally $\mathfrak{P} \notin T'$. If $z_\mathfrak{P}$ is unramified, then it is cyclic. Let $z_\mathfrak{P}$ be ramified. Since $z_\mathfrak{P} = z'_\mathfrak{P} + z''_\mathfrak{P}$ and $z'_\mathfrak{P}$ is unramified, $z''_\mathfrak{P}$ must be ramified. Thus $\mathfrak{p}$ splits completely in $\Omega'$, hence in $K'$. From this we obtain $z'_{\sigma\mathfrak{P}} = 0$ for all $\sigma \in G(K|k)$, since by definition of $K'$ the element $z'$ becomes zero in $H^1(K', A')$ and thus $z'_{\sigma\mathfrak{P}}$ is zero in $H^1(K'_\mathfrak{Q}, A') = H^1(K_{\sigma\mathfrak{P}}, A')$, where $\mathfrak{Q}$ is a prime of $K'$ above $\sigma\mathfrak{P}$. Therefore $z_\mathfrak{P} = z''_\mathfrak{P}$, i.e. $z_\mathfrak{P}$ is cyclic and $z_{\sigma\mathfrak{P}} = z''_{\sigma\mathfrak{P}} = 0$ for $\sigma \in G(K|k)\setminus\{1\}$.

It therefore remains to construct a class $\tilde{y}$ with the above properties. Consider the commutative and exact diagram

$$\begin{array}{ccc}
\prod_{S'\setminus T'} H^1(K_\mathfrak{P}, \mu_p) & \xrightarrow{\sim} & \Big(\prod_{S'\setminus T'} H^1(K_\mathfrak{P}, \mathbb{Z}/p\mathbb{Z})\Big)^\vee \\
\uparrow & & \uparrow \\
H^1(k_{S'}|K, \mu_p) \longrightarrow \prod_{S'} H^1(K_\mathfrak{P}, \mu_p) & \longrightarrow & H^1(k_{S'}|K, \mathbb{Z}/p\mathbb{Z})^\vee \\
\downarrow & & \downarrow \\
\prod_{T'} H^1(K_\mathfrak{P}, \mu_p) & \xrightarrow{\alpha} & H^1(\Omega'|K, \mathbb{Z}/p\mathbb{Z})^\vee.
\end{array}$$

If we can show that $\alpha$ annihilates the element $\xi = (\xi_\mathfrak{P})_{\mathfrak{P} \in T'(K)} \in \prod_{T'} H^1(K_\mathfrak{P}, \mu_p)$ given by $\xi_\mathfrak{P} = y''_\mathfrak{P}$ for $\mathfrak{P} \in T(K)$ and $\xi_\mathfrak{P} = 0$ for $\mathfrak{P} \in T' \setminus T(K)$, then the existence of $\tilde{y}$ follows by diagram chasing. We use the injection

$$H^1(\Omega'|K, \mathbb{Z}/p\mathbb{Z})^\vee \hookrightarrow H^1(\Omega|K, \mathbb{Z}/p\mathbb{Z})^\vee \oplus H^1(\tilde{K}|K, \mathbb{Z}/p\mathbb{Z})^\vee$$

in order to write the image of $\xi$ in the form $\alpha(\xi) = (\alpha_1(\xi), \alpha_2(\xi))$. Since $\tilde{K}|K$ is unramified at all $\mathfrak{P} \in T(K)$, $\alpha_2$ factorizes over the quotient

$$\prod_{T'} H^1(K_\mathfrak{P}, \mu_p) / \Big(\prod_T H^1_{nr}(K_\mathfrak{P}, \mu_p) \times \prod_{T'\setminus T} \{0\}\Big).$$

Hence $\alpha_2(\xi) = 0$. Finally the diagram

$$\begin{array}{ccc}
\prod_{T'} H^1(K_\mathfrak{P}, \mu_p) & \xrightarrow{\alpha} & H^1(\Omega'|K, \mathbb{Z}/p\mathbb{Z})^\vee \\
\downarrow & & \downarrow \\
\prod_T H^1(K_\mathfrak{P}, \mu_p) & \longrightarrow & H^1(\Omega|K, \mathbb{Z}/p\mathbb{Z})^\vee \\
\uparrow & & \uparrow \\
H^1(k_S|K, \mu_p) \longrightarrow \prod_S H^1(K_\mathfrak{P}, \mu_p) & \longrightarrow & H^1(k_S|K, \mathbb{Z}/p\mathbb{Z})^\vee
\end{array}$$

shows that $\alpha_1(\xi)$ is equal to the image of $y''$ ($\in H^1(k_S|K, \mu_p)$) in $H^1(\Omega|K, \mathbb{Z}/p\mathbb{Z})^\vee$, hence trivial by the exactness of the lower line. This finishes the proof. □

## 5. Proof of Šafarevic's theorem

The essential step in the proof of theorem 1 is the following

**Theorem 14** *Let $K|k$ be a finite Galois extension of the global field $k$ and let $\varphi : G_k \twoheadrightarrow G(K|k) = G$. Then every split embedding problem*

$$
\begin{array}{c}
G_k \\
\downarrow \\
1 \longrightarrow H \longrightarrow H \rtimes G \rightleftarrows G \longrightarrow 1
\end{array}
$$

*with finite nilpotent kernel $H$ has a proper solution.*

Since a finite nilpotent group is the direct product of its $p$-Sylow subgroups and since every finite $G$-operator $p$-group is a quotient of $\mathcal{F}(n)/\mathcal{F}(n)^{(\nu)}$ for some $n, \nu$, it suffices to show the following assertion:

*For every prime number $p$, all $n \in \mathbb{N}$ and all $\nu = (i, j)$ the split embedding problem*

$$
\begin{array}{c}
G_k \\
\downarrow \\
1 \longrightarrow \mathcal{F}(n)/\mathcal{F}(n)^{(\nu)} \longrightarrow \mathcal{F}(n)/\mathcal{F}(n)^{(\nu)} \rtimes G \rightleftarrows G \longrightarrow 1
\end{array}
$$

*has a proper solution $N_\nu^n|k$.*

Let us first assume that $p \neq \mathrm{char}(k)$. We will proceed by induction on $\nu$ whereas $n$ will be arbitrary. If $\nu = (1,1)$ nothing is to show. Now we assume that we have already found a solution $\varphi_{n,\nu} : G_k \twoheadrightarrow \mathcal{F}(n)/\mathcal{F}(n)^{(\nu)} \rtimes G$ and we consider the embedding problem

$$(*) \qquad\qquad\qquad\qquad\qquad\qquad\qquad\qquad G_k$$
$$\downarrow \varphi_{n,\nu}$$
$$\mathcal{F}(n)^{(\nu)}/\mathcal{F}(n)^{(\nu+1)} \hookrightarrow \mathcal{F}(n)/\mathcal{F}(n)^{(\nu+1)} \rtimes G \twoheadrightarrow \mathcal{F}(n)/\mathcal{F}(n)^{(\nu)} \rtimes G.$$

This embedding problem is in general not solvable, but we are going to solve it after replacing $\varphi_{n,\nu}$ by another solution $\tilde{\varphi}_{n,\nu}$ for the induction step $\nu$ on the level $n$. This new solution is induced by a solution

$$\varphi_{m,\nu} : G_k \to \mathcal{F}(m)/\mathcal{F}(m)^{(\nu)} \rtimes G$$

for some large $m \geq n$ via a suitable chosen $G$-invariant surjection $\psi : \mathcal{F}(m) \twoheadrightarrow \mathcal{F}(n)$. Let us consider the associated commutative and exact diagram

$$(**) \qquad \begin{array}{ccccc}
 & & & & G_k \\
 & & & & \downarrow \varphi_{m,\nu} \\
\mathcal{F}(m)^{(\nu)}/\mathcal{F}(m)^{(\nu+1)} & \hookrightarrow & \mathcal{F}(m)/\mathcal{F}(m)^{(\nu+1)} \rtimes G & \twoheadrightarrow & \mathcal{F}(m)/\mathcal{F}(m)^{(\nu)} \rtimes G \\
\downarrow \bar{\psi} & & \downarrow \psi_{\nu+1} & & \downarrow \psi_\nu \\
\mathcal{F}(n)^{(\nu)}/\mathcal{F}(n)^{(\nu+1)} & \hookrightarrow & \mathcal{F}(n)/\mathcal{F}(n)^{(\nu+1)} \rtimes G & \twoheadrightarrow & \mathcal{F}(n)/\mathcal{F}(n)^{(\nu)} \rtimes G.
\end{array}$$

To shorten notations we set again:

$$\mathcal{F}(n)/\nu = \mathcal{F}(n)/\mathcal{F}(n)^{(\nu)}, \quad \mathcal{E}(n,\nu) = \mathcal{F}(n)^{(\nu)}/\mathcal{F}(n)^{(\nu+1)}.$$

Since $\mathcal{E}(n,\nu)$ is contained in the center of $\mathcal{F}(n)/\nu+1$, the action of $\mathcal{F}(n)/\nu \rtimes G$ on $\mathcal{E}(n,\nu)$ factors over the canonical projection

$$\mathcal{F}(n)/\nu \rtimes G \twoheadrightarrow G,$$

in particular $G_K \subseteq G_k$ acts trivially on $\mathcal{E}(n,\nu)$.

Let $\alpha_m$ and $\alpha_n$ denote the 2-classes corresponding to the group extensions in $(**)$ and consider the commutative exact diagram

$$\begin{array}{ccc}
H^2(\mathcal{F}(m)/\nu \rtimes G, \mathcal{E}(m,\nu)) & \xrightarrow{\varphi_{m,\nu}^*} & H^2(k, \mathcal{E}(m,\nu)) \\
\bar{\psi}_* \downarrow & & \downarrow \bar{\psi}_* \\
H^2(\mathcal{F}(m)/\nu \rtimes G, \mathcal{E}(n,\nu)) & \xrightarrow{\varphi_{m,\nu}^*} & H^2(k, \mathcal{E}(n,\nu)) \\
\psi_\nu^* \uparrow & & \| \\
H^2(\mathcal{F}(n)/\nu \rtimes G, \mathcal{E}(n,\nu)) & \xrightarrow{\varphi_{n,\nu}^*} & H^2(k, \mathcal{E}(n,\nu)).
\end{array}$$

Now [1] chap.13 th.2 asserts that $\psi_\nu^*(\alpha_n) = \bar{\psi}_*(\alpha_m)$ and by [3] Satz 1.1 the embedding problem on the level $n$ is solvable if and only if $\varphi_{n,\nu}^*(\alpha_n) = 0$. We are searching for an $m \geq n$, a solution $\varphi_{m,\nu}$ on the level $m$ and a suitable surjective $G$-homomorphism $\psi : \mathcal{F}(m) \twoheadrightarrow \mathcal{F}(n)$, such that

$$\bar{\psi}_*(\varphi_{m,\nu}^*(\alpha_m)) = 0 \in H^2(k, \mathcal{E}(n,\nu)).$$

As we will show below, this can be achieved for large enough $m$ if the solution $\varphi_{m,\nu}$ is of a special type. If we can guarantee that the new solution is also of this special type, then the induction process works for a modified, stronger statement. In fact we are going to prove the sharpened

**Theorem 15** *Let $K|k$ be a finite Galois extension of the global field $k$ and let $\varphi : G_k \twoheadrightarrow G(K|k) = G$. Then for every prime number $p$, all $n \in \mathbb{N}$ and all $\nu = (i, j)$ the split embedding problem*

$$
\begin{array}{c}
G_k \\
\downarrow \varphi \\
1 \longrightarrow \mathcal{F}(n)/\mathcal{F}(n)^{(\nu)} \longrightarrow \mathcal{F}(n)/\mathcal{F}(n)^{(\nu)} \rtimes G \rightleftarrows G \longrightarrow 1
\end{array}
$$

*has a proper solution $N_\nu^n | k$. If $p \neq \operatorname{char} k$, we can choose the solution in such a way that the following conditions are satisfied:*

(i) *All $\mathfrak{p} \in \operatorname{Ram}(K|k) \cup S_p \cup S_\infty$ are completely decomposed in $N_\nu^n | K$.*

(ii) *If $\mathfrak{p}$ is ramified in $N_\nu^n | K$, then $\mathfrak{p}$ splits completely in $K|k$ and $N_{\nu,\mathfrak{p}}^n | k_\mathfrak{p}$ is a (cyclic) totally ramified extension of local fields.*

**Proof:** We prove the theorem by induction on $\nu$, whereas $n$ and $G$ are arbitrary. We defer the case $\operatorname{char}(k) = p$ and assume that $\operatorname{char}(k) \neq p$. If $\nu = (1, 1)$, nothing is to show. We prove the induction step, i.e. we solve the embedding problem defined by the diagram (∗) above in four substeps. Furthermore, for the induction step $v \mapsto v + 1$, we may assume that $\mu_{p^e} \subseteq K$, where $p^e$ is the exponent of the group $\mathcal{F}(n)/\mathcal{F}(n)^{(\nu+1)}$ (which does depend on $\nu$ but not on $n$ and $G$). Otherwise we lift the embedding problem via $G_k \twoheadrightarrow G(K(\mu_{p^e})|k) \twoheadrightarrow G(K|k)$, thus making $G$ bigger. Note that this does not influence conditions (i) and (ii).

**First Step:** *The problem (∗) induces local split embedding problems at all $\mathfrak{p} \in \operatorname{Ram}(K|k) \cup S_p \cup S_\infty$ and is locally solvable (not necessarily proper) at every prime $\mathfrak{p}$ after changing $\varphi_{n,\nu}$.*

a) If $\mathfrak{p} \in \operatorname{Ram}(K|k) \cup S_p \cup S_\infty$, then $G_\mathfrak{p}(N_\nu^n | k) = (\mathcal{F}(n)/\mathcal{F}(n)^{(\nu)} \rtimes G)_\mathfrak{p} \cong G_\mathfrak{p}(K|k)$ by (i). We show that, after changing $\varphi_{n,\nu}$, the local group extensions corresponding to these primes are split extensions. In particular, the associated local embedding problems are solvable in a trivial way.

Let $\alpha_n(\mathfrak{p})$ be the 2-class in $H^2((\mathcal{F}(n)/\mathcal{F}(n)^{(\nu)} \rtimes G)_\mathfrak{p}, \mathcal{F}(n)^{(\nu)}/\mathcal{F}(n)^{(\nu+1)})$ which corresponds to the group extension given by the upper line of the diagram

$$
\begin{array}{ccccc}
\mathcal{F}(n)^{(\nu)}/\mathcal{F}(n)^{(\nu+1)} & \hookrightarrow & E_\mathfrak{p} & \twoheadrightarrow & (\mathcal{F}(n)/\mathcal{F}(n)^{(\nu)} \rtimes G)_\mathfrak{p} \\
\| & & \downarrow & & \downarrow \\
\mathcal{F}(n)^{(\nu)}/\mathcal{F}(n)^{(\nu+1)} & \hookrightarrow & \mathcal{F}(n)/\mathcal{F}(n)^{(\nu+1)} \rtimes G & \twoheadrightarrow & \mathcal{F}(n)/\mathcal{F}(n)^{(\nu)} \rtimes G \,.
\end{array}
$$

Apply the induction hypothesis to the corresponding embedding problem on a large level $m$. The number $m$ and a surjective $G$-invariant homomorphism $\psi : \mathcal{F}(m) \twoheadrightarrow \mathcal{F}(n)$ will be chosen below. By (i) we know that $\mathfrak{p}$ is completely decomposed in $N_\nu^m | K$ and we have a commutative diagram for the associated local groups (writing $G_\mathfrak{p}$ for $G_\mathfrak{p}(K|k)$):

$$\begin{array}{ccc} (\mathcal{F}(m)/\nu \rtimes G)_\mathfrak{p} & \xrightarrow{\sim} & G_\mathfrak{p} \\ \downarrow \wr & & \| \\ (\mathcal{F}(n)/\nu \rtimes G)_\mathfrak{p} & \xrightarrow{\sim} & G_\mathfrak{p} \end{array}$$

Therefore we obtain the diagram

$$\begin{array}{ccc} H^2((\mathcal{F}(m)/\nu \rtimes G)_\mathfrak{p}, \mathcal{E}(m,\nu)) & \xleftarrow[\sim]{\text{inf}} & H^2(G_\mathfrak{p}, \mathcal{E}(m,\nu)) \\ \downarrow \bar\psi_* & & \downarrow \bar\psi_* \\ H^2((\mathcal{F}(m)/\nu \rtimes G)_\mathfrak{p}, \mathcal{E}(n,\nu)) & \xleftarrow[\sim]{\text{inf}} & H^2(G_\mathfrak{p}, \mathcal{E}(n,\nu)) \\ \updownarrow \psi_\nu^* & & \| \\ H^2((\mathcal{F}(n)/\nu \rtimes G)_\mathfrak{p}, \mathcal{E}(n,\nu)) & \xleftarrow[\sim]{\text{inf}} & H^2(G_\mathfrak{p}, \mathcal{E}(n,\nu)) \end{array}$$

Using proposition 6 with $G$, $k = 2$ and $T = \mathrm{Ind}_G^{G_\mathfrak{p}} \mathbb{F}_p$ we find an $m \geq n$ and a surjective pro-$p$-$G$ homomorphism $\psi : \mathcal{F}(m) \twoheadrightarrow \mathcal{F}(n)$ such that the homomorphism

$$\begin{array}{rcl} H^2(G_\mathfrak{p}, \mathcal{E}(m,\nu)) & = & H^2(G, \mathcal{E}(m,\nu) \otimes T)) \\ & & \downarrow \bar\psi_* \\ H^2(G_\mathfrak{p}, \mathcal{E}(n,\nu)) & = & H^2(G, \mathcal{E}(n,\nu) \otimes T) \end{array}$$

maps $\mathit{inf}^{-1}(\alpha_m(\mathfrak{p}))$ to 0, hence the above diagram implies that $\alpha_n(\mathfrak{p}) = 0$.

Now we can execute the above procedure for all the finitely many primes $\mathfrak{p} \in \mathrm{Ram}(K|k) \cup S_p \cup S_\infty$, making $m$ each times bigger. Note that we do not destroy the already achieved success for the primes say $\mathfrak{p}_1, \ldots, \mathfrak{p}_r$ for which the local embedding problems already split. Indeed, the property of inducing a split embedding problem at a prime $\mathfrak{p}$ survives the shrinking process from $m$ to $n$ if $\mathfrak{p}$ is completely decomposed is $N_\nu^m | K$ (and we suppose this for $\mathfrak{p} \in \mathrm{Ram}(K|k) \cup S_p \cup S_\infty$). Therefore we can execute a shrinking process for $\mathfrak{p}_{r+1}$, inducing a solution from $\mathcal{F}(m)$, which has already the desired property for $\mathfrak{p}_1, \ldots, \mathfrak{p}_r$ (and which we could have produced by another shrinking $\mathcal{F}(m') \twoheadrightarrow \mathcal{F}(m)$).

An alternative way to proceed at this point is to replace the module $T$ in the above argument by the direct sum of $\mathrm{Ind}_G^{G_\mathfrak{p}} \mathbb{F}_p$, where $\mathfrak{p}$ runs through $\mathrm{Ram}(K|k) \cup S_p \cup S_\infty$. In this way we deal with all these primes within one shrinking process.

Finally we see that the embedding problem

$$
\begin{array}{c}
G_k \\
\downarrow \varphi_{m,\nu} \\
\mathcal{F}(m)/\mathcal{F}(m)^{(\nu)} \rtimes G \\
\downarrow \bar{\psi}
\end{array}
$$

$$\mathcal{F}(n)^{(\nu)}/\mathcal{F}(n)^{(\nu+1)} \hookrightarrow \mathcal{F}(n)/\mathcal{F}(n)^{(\nu+1)} \rtimes G \longrightarrow \mathcal{F}(n)/\mathcal{F}(n)^{(\nu)} \rtimes G$$

induces split group extensions for all $\mathfrak{p} \in \mathrm{Ram}(K|k) \cup S_p \cup S_\infty$.

*We have just changed $\varphi_{n,\nu}$ to $\bar{\psi} \circ \varphi_{m,\nu}$, thus $N_\nu^n$ to some $\tilde{N}_\nu^n$. But by assumption conditions (i) and (ii) are also satisfied for $N_\nu^m$ hence for the new field $\tilde{N}_\nu^n$. We will not change the notations.*

b) If $\mathfrak{p}$ is unramified in $N_\nu^n|k$, then the homomorphism

$$\hat{\mathbb{Z}} \cong G_\mathfrak{p}(k)/T_\mathfrak{p}(k) \xrightarrow{\varphi_\mathfrak{p}} G_\mathfrak{p}(N_\nu^n|k)$$

extends obviously since $\hat{\mathbb{Z}}$ is free.

c) Let $\mathfrak{p} \in \mathrm{Ram}(N_\nu^n|K)$. Then $\mathfrak{p}$ splits completely in $K|k$ and $G_\mathfrak{p}(N_\nu^n|K) \cong \mathbb{Z}/p^a\mathbb{Z}$ by condition (ii). We assumed that $\mu_{p^e} \subseteq K$, where $p^e$ is the exponent of the group $\mathcal{F}(n)/\mathcal{F}(n)^{(\nu+1)}$. Since $N_{\nu,\mathfrak{p}}^n|K_\mathfrak{p}$ is totally ramified by assumption, there exists a prime element $\pi_\mathfrak{p}$ of $K_\mathfrak{p}$ such that $N_{\nu,\mathfrak{p}}^n = K_\mathfrak{p}(\sqrt[p^a]{\pi_\mathfrak{p}})$. An arbitrary chosen preimage of a generator of the cyclic group $G_\mathfrak{p}(N_\nu^n|K)$ in $\mathcal{F}(n)/\mathcal{F}(n)^{(\nu+1)} \rtimes G$ has order $p^{a+\varepsilon}$, where $0 \leq \varepsilon \leq 1$. We can solve our embedding problem by taking a $p^{a+\varepsilon}$-th root of $\pi_\mathfrak{p}$, since $\mu_{p^{a+\varepsilon}} \subseteq \mu_{p^e} \subseteq K \subseteq K_\mathfrak{p}$.

**Second Step:** *The problem $(*)$ induces local split embedding problems at all $\mathfrak{p} \in \mathrm{Ram}(K|k) \cup S_p \cup S_\infty$ and is globally solvable (not necessarily proper) after changing $\varphi_{n,\nu}$.*

As above we consider the problem for different numbers $m \geq n$:

$$
\begin{array}{ccccc}
& & & & G_k \\
& & & & \downarrow \varphi_{m,\nu} \\
\mathcal{F}(m)^{(\nu)}/\mathcal{F}(m)^{(\nu+1)} & \hookrightarrow & \mathcal{F}(m)/\mathcal{F}(m)^{(\nu+1)} \rtimes G & \longrightarrow & \mathcal{F}(m)/\mathcal{F}(m)^{(\nu)} \rtimes G \\
\downarrow \bar{\psi} & & \downarrow \psi_{\nu+1} & & \downarrow \psi_\nu \\
\mathcal{F}(n)^{(\nu)}/\mathcal{F}(n)^{(\nu+1)} & \hookrightarrow & \mathcal{F}(n)/\mathcal{F}(n)^{(\nu+1)} \rtimes G & \longrightarrow & \mathcal{F}(n)/\mathcal{F}(n)^{(\nu)} \rtimes G.
\end{array}
$$

Let $\alpha_m$ and $\alpha_n$ denote the 2-classes corresponding to the above group extensions. A surjective pro-$p$-$G$ homomorphism $\psi : \mathcal{F}(m) \twoheadrightarrow \mathcal{F}(n)$ (inducing $\psi_\nu$, $\psi_{\nu+1}$ and $\bar{\psi}$) will be defined below. Both problems (for $m$ with $\varphi_{m,\nu}$ and for $n$ with $\varphi_{n,\nu} = \psi_\nu \circ \varphi_{m,\nu}$) are locally solvable by step 1.

In order to show the existence of a solution $\varphi_{n,\nu+1}$ we have to prove that the 2-class $\alpha_n$ maps to zero under the inflation map $\varphi_{n,\nu}^*$, see [3] Satz 1.1:

$$\begin{array}{ccc}
 & & \prod_{\mathfrak{p}} H^2(k_{\mathfrak{p}}, \mathcal{E}(n,\nu)) \\
 & & \uparrow \\
H^2(\mathcal{F}(n)/\mathcal{F}(n)^{(\nu)} \rtimes G, \mathcal{E}(n,\nu)) & \xrightarrow{\varphi_{n,\nu}^*} & H^2(k, \mathcal{E}(n,\nu)) \\
 & & \uparrow \\
 & & \text{III}^2(k, \mathcal{E}(n,\nu)).
\end{array}$$

Here we have set as before $\mathcal{E}(n,\nu) = \mathcal{F}(n)^{(\nu)}/\mathcal{F}(n)^{(\nu+1)}$. By the first step we can assume that

$$\varphi_{n,\nu}^*(\alpha_n) \in \text{III}^2(k, \mathcal{E}(n,\nu)) \quad , \quad \varphi_{m,\nu}^*(\alpha_m) \in \text{III}^2(k, \mathcal{E}(m,\nu)).$$

As we have already observed at the beginning

$$\psi_\nu^*(\alpha_n) = \bar{\psi}_*(\alpha_m)$$

which gives us

$$\begin{array}{rcl}
\varphi_{n,\nu}^*(\alpha_n) & = & (\psi_\nu \circ \varphi_{m,\nu})^*(\alpha_n) = (\varphi_{m,\nu}^* \circ \psi_\nu^*)(\alpha_n) \\
 & = & (\varphi_{m,\nu}^* \circ \bar{\psi}_*)(\alpha_m) = (\bar{\psi}_* \circ \varphi_{m,\nu}^*)(\alpha_m).
\end{array}$$

In order to shrink the obstruction, we look for a surjective homomorphism onto $\text{III}^2(k, \mathcal{E}(n,\nu))$, which has a "shrinkable" source, i.e. a source to which proposition 6 applies.

*Claim:* We have a commutative diagram

$$\begin{array}{ccc}
H^{-2}(G, \mathcal{E}(m,\nu)(-1)) & \twoheadrightarrow & \text{III}^2(k, \mathcal{E}(m,\nu)) \\
\downarrow \bar{\psi}_* & & \downarrow \bar{\psi}_* \\
H^{-2}(G, \mathcal{E}(n,\nu)(-1)) & \twoheadrightarrow & \text{III}^2(k, \mathcal{E}(n,\nu))
\end{array}$$

with surjective horizontal maps, where $(-1)$ denotes the $(-1)$-Tate twist.

Using this claim and proposition 6 with $k = -2$, the $\mathbb{F}_p[G]$-module $T = \text{Hom}(\mu_p, \mathbb{Z}/p\mathbb{Z})$ and an element $x_1$ which is a pre-image of $\varphi_{m,\nu}^*(\alpha_m)$ in the group $H^{-2}(G, \mathcal{E}(m,\nu)(-1))$, we obtain a surjective pro-$p$-$G$ operator homomorphism $\psi : \mathcal{F}(m) \twoheadrightarrow \mathcal{F}(n)$ such that $\bar{\psi}_*\varphi_{m,\nu}^*(\alpha_m) = 0$, hence $\varphi_{n,\nu}^*(\alpha_n) = 0$. Thus the embedding problem is solvable. Furthermore, as explained in the first step, the local condition at the primes in $\text{Ram}(K|k) \cup S_p \cup S_\infty$ remains untouched within the shrinking process. In order to finish step 2 it remains to give the

*Proof of the claim:* Let $\mathcal{E}(n,\nu)' = \mathrm{Hom}(\mathcal{E}(n,\nu), \mu_p)$. By the Tate-Poitou duality theorem we know that

$$\mathrm{III}^2(k, \mathcal{E}(n,\nu)) \cong \mathrm{III}^1(k, \mathcal{E}(n,\nu)')^* \,.$$

Using the Hasse-principle and the fact that $\mathcal{E}(n,\nu)'$ is a trivial $G_K$-module ($\mu_p \subseteq K$) we obtain that the homomorphisms $\iota$ in the commutative and exact diagram

$$\begin{array}{ccc}
H^1(K, \mathcal{E}(n,\nu)') & \stackrel{\iota}{\hookrightarrow} & \prod_{\mathfrak{P}} H^1(K_{\mathfrak{P}}, \mathcal{E}(n,\nu)') \\
\uparrow & & \uparrow \\
0 \to \mathrm{III}^1(k, \mathcal{E}(n,\nu)') \longrightarrow H^1(k, \mathcal{E}(n,\nu)') & \longrightarrow & \prod_{\mathfrak{p}} H^1(k_{\mathfrak{p}}, \mathcal{E}(n,\nu)') \\
\uparrow & & \\
H^1(K|k, \mathcal{E}(n,\nu)') & & \\
\uparrow & & \\
0 & &
\end{array}$$

is injective. Hence we get an injection $\mathrm{III}^1(k, \mathcal{E}(n,\nu)') \hookrightarrow H^1(K|k, \mathcal{E}(n,\nu)')$ and regarding that the dual of cohomology is homology (hence for finite groups cohomology in negative dimensions), we obtain a canonical surjection

$$\begin{array}{ccc}
H^1(G, \mathcal{E}(n,\nu)')^* & \twoheadrightarrow & \mathrm{III}^1(k, \mathcal{E}(n,\nu)')^* \\
\| & & \| \\
H^{-2}(G, \mathcal{E}(n,\nu)(-1)) & & \mathrm{III}^2(k, \mathcal{E}(n,\nu)) \,.
\end{array}$$

This proves the claim.

**Third step:** *After changing $\varphi_{n,\nu}$, the problem $(*)$ has a proper global solution, which satisfies condition (i) and all primes $\mathfrak{p} \in \mathrm{Ram}(N^n_{\nu+1}|N^n_\nu) \smallsetminus \mathrm{Ram}(N^n_\nu|K)$ are completely decomposed in $N^n_\nu|k$. Furthermore, the local extension $N^n_{\nu+1,\mathfrak{p}}|K_{\mathfrak{p}}$ is (cyclic) totally ramified for $\mathfrak{p} \in \mathrm{Ram}(N^n_\nu|K)$.*

We achieve this by the following procedure: Consider a solution $\varphi_{n,\nu+1}$ of the embedding problem $(*)$ which we have obtained in step 2. Its equivalence class $[\varphi_{n,\nu+1}]$ is an element of the space $\mathscr{S}_{(*)}$ of solutions of $(*)$ modulo equivalence. Conditions (i),(ii), properness and the other conditions that we want to achieve in this third step only depend on the equivalence class of a solution. The space $\mathscr{S}_{(*)}$ is a principal homogeneous space over $H^1(G_k, \mathcal{E}(n,\nu))$ (see [6] Satz 1.3). Recall that the action is defined as follows: Choose a representing cocycle $G_k \to \mathcal{E}(n,\nu)$ and multiply a solution $G_k \to \mathcal{F}(n)/\nu+1 \rtimes G$ of the embedding problem with the cocycle. This yields a map $G_k \to \mathcal{F}(n)/\nu+1 \rtimes G$ which is a homomorphism (!) and the equivalence class of this new solution is independent of the made choices.

Now we look for a suitable cohomology class $\varepsilon \in H^1(G_k, \mathcal{E}(n,\nu))$, such that the new solution
$$\tilde{\varphi}_{n,\nu+1} = {}^{\varepsilon}\varphi_{n,\nu+1}$$
has the required properties. We assume that $\varphi_{m,\nu+1}$ is obtained from step 2 and that $\varphi_{m,\nu}$ satisfies (i) and (ii). Note that the properness of the solution is only a problem for the first step $(1,1) \to (2,1)$ since in all higher induction steps the properness follows automatically from the induction hypothesis and from the Frattini argument.

Let us consider how the local behaviour of $\tilde{\varphi}_{n,\nu+1} = {}^{\varepsilon}\varphi_{n,\nu+1}$ is connected to that of $\varphi_{n,\nu+1}$. By this we mean, that we want to compare the ramification and decomposition of primes in the associated field extensions $\tilde{N}^n_{\nu+1}|N^n_\nu$ and $N^n_{\nu+1}|N^n_\nu$. (Since we do not know, whether the solutions are proper, one or both of these field extensions might be trivial.) Let $\mathfrak{p}$ be a prime in $N^n_\nu$. The behaviour of $\mathfrak{p}$ in $N^n_{\nu+1}|N^n_\nu$ is characterized by the homomorphism
$$\varphi_{n,\nu+1}|_{G_{(N^n_\nu)_\mathfrak{p}}} \in \mathrm{Hom}(G_{(N^n_\nu)_\mathfrak{p}}, \mathcal{E}(n,\nu)).$$

Since $G_{(N^n_\nu)_\mathfrak{p}}$ acts trivially on $\mathcal{E}(n,\nu)$, we can interpret $\varphi_{n,\nu+1}|_{G_{(N^n_\nu)_\mathfrak{p}}}$ as an element in
$$H^1((N^n_\nu)_\mathfrak{p}, \mathcal{E}(n,\nu))^{G_{k_\mathfrak{p}}}.$$

Consider the exact sequence
$$0 \longrightarrow H^1((N^n_\nu)_\mathfrak{p}|k_\mathfrak{p}) \longrightarrow H^1(k_\mathfrak{p}) \xrightarrow{\alpha} H^1((N^n_\nu)_\mathfrak{p})^{G_{k_\mathfrak{p}}} \xrightarrow{\beta} H^2((N^n_\nu)_\mathfrak{p}|k_\mathfrak{p}),$$
which is obtained from the Hochschild-Serre sequence for the tower $\bar{k}_\mathfrak{p} \mid (N^n_\nu)_\mathfrak{p} \mid k_\mathfrak{p}$ and in which $\mathcal{E}(n,\nu)$ are the (not written) coefficients of the cohomology groups. We see that $\tilde{\varphi}_{n,\nu+1}|_{G_{(N^n_\nu)_\mathfrak{p}}}$ is given by
$$\varphi_{n,\nu+1}|_{G_{(N^n_\nu)_\mathfrak{p}}} + \alpha(\varepsilon) \in H^1((N^n_\nu)_\mathfrak{p}, \mathcal{E}(n,\nu))^{G_{k_\mathfrak{p}}}.$$

Now we choose a finite set $T^0$ of primes in $cs(N^n_\nu|k)$ and homomorphisms
$$x_\mathfrak{p} : G_{k_\mathfrak{p}}/T_{k_\mathfrak{p}} = G_{(N^n_\nu)_\mathfrak{p}}/T_{(N^n_\nu)_\mathfrak{p}} \longrightarrow \mathcal{E}(n,\nu)$$
for $\mathfrak{p} \in T^0$ such that their images generate $\mathcal{E}(n,\nu)$. (The set $T^0$ will be responsible for the properness of the new solution.)

Set
$$\begin{aligned}
T^1 &= \mathrm{Ram}(K|k) \cup S_p \cup S_\infty, \\
T^2 &= \mathrm{Ram}(N^n_\nu|K), \\
T^3 &= \mathrm{Ram}(N^n_{\nu+1}|K)\backslash(\mathrm{Ram}(N^n_\nu|k) \cup S_p \cup S_\infty), \\
T &= T^0 \cup T^1 \cup T^2 \cup T^3 \quad \text{and} \\
S &= cs(N^n_\nu|k) \cup T.
\end{aligned}$$

Since $\mathfrak{p} \in T^0$ splits completely in $N_\nu^n | k$, there exists a $\xi_\mathfrak{p} \in H^1(k_\mathfrak{p}, \mathcal{E}(n,\nu))$ with

$$\alpha(\xi_\mathfrak{p}) = x_\mathfrak{p} : G_{(N_\nu^n)_\mathfrak{p}} \longrightarrow \mathcal{E}(n,\nu).$$

Let $\mathfrak{p} \in T^1$. Then, by the induction hypothesis, the extension $(N_\nu^n)_\mathfrak{p} | K_\mathfrak{p}$ is trivial and the group extension in the diagram

$$\begin{array}{ccccccccc}
 & & & & & & G_{k_\mathfrak{p}} & & \\
 & & & \varphi_{n,\nu+1}|_{G_{k_\mathfrak{p}}} & \swarrow & & \downarrow \varphi_{n,\nu}|_{G_{k_\mathfrak{p}}} & & \\
1 & \longrightarrow & \mathcal{E}(n,\nu) & \longrightarrow & E_\mathfrak{p} & \longrightarrow & G_\mathfrak{p}(N_\nu^n | k) & \longrightarrow & 1
\end{array}$$

splits. Hence $\beta(\varphi_{n,\nu+1}|_{G_{(N_\nu^n)_\mathfrak{p}}}) = 0$ and we therefore find $\xi_\mathfrak{p} \in H^1(k_\mathfrak{p}, \mathcal{E}(n,\nu))$ such that

$$\alpha(\xi_\mathfrak{p}) = -\varphi_{n,\nu+1}|_{G_{(N_\nu^n)_\mathfrak{p}}} : G_{(N_\nu^n)_\mathfrak{p}} \longrightarrow \mathcal{E}(n,\nu).$$

If $\mathfrak{p} \in T^2$, then by the induction hypothesis $\mathfrak{p} \notin T^1$, $K_\mathfrak{p} | k_\mathfrak{p}$ is trivial and $(N_\nu^n)_\mathfrak{p} | K_\mathfrak{p}$ is a (cyclic) totally ramified extension. We consider the exact and commutative diagram,

$$\begin{array}{ccc}
H^1_{nr}(k_\mathfrak{p}) & \dashrightarrow & H^1_{nr}((N_\nu^n)_\mathfrak{p})^{G_{k_\mathfrak{p}}} \\
\downarrow & & \downarrow \\
H^1((N_\nu^n)_\mathfrak{p} | k_\mathfrak{p}) \hookrightarrow H^1(k_\mathfrak{p}) & \xrightarrow{\alpha} & H^1((N_\nu^n)_\mathfrak{p})^{G_{k_\mathfrak{p}}}.
\end{array}$$

Since $\mathfrak{p} \in T^2$, the dotted arrow in the diagram above is an isomorphism. Thus there is a $\xi_\mathfrak{p} \in H^1(k_\mathfrak{p}, \mathcal{E}(n,\nu))$ such that

$$\varphi_{n,\nu+1}|_{G_{(N_\nu^n)_\mathfrak{p}}} + \alpha(\xi_\mathfrak{p}) : G_{(N_\nu^n)_\mathfrak{p}} \longrightarrow \mathcal{E}(n,\nu)$$

is either trivial (if $\varphi_{n,\nu+1}|_{G_{(N_\nu^n)_\mathfrak{p}}}$ is unramified) or induces a totally ramified extension of degree $p$ of $(N_\nu^n)_\mathfrak{p}$.

If $\mathfrak{p} \in T^3$, then we have a commutative and exact diagram

$$\begin{array}{ccccc}
H^1((N_\nu^n)_\mathfrak{p} | k_\mathfrak{p}) & \hookrightarrow & H^1(k_\mathfrak{p}) & \xrightarrow{\alpha} & H^1((N_\nu^n)_\mathfrak{p})^{G_{k_\mathfrak{p}}} \\
& & \downarrow & & \downarrow {\scriptstyle res} \\
& & H^1(T_{k_\mathfrak{p}})^{G_{k_\mathfrak{p}}} & \dashrightarrow & H^1(T_{(N_\nu^n)_\mathfrak{p}})^{G_{k_\mathfrak{p}}}
\end{array}$$

where now the lower dotted arrow is an isomorphism, since $(N_\nu^n)_\mathfrak{p} | k_\mathfrak{p}$ is unramified. Let $\xi_\mathfrak{p} \in H^1(k_\mathfrak{p}, \mathcal{E}(n,\nu))$ such that

$$\varphi_{n,\nu+1}|_{G_{(N_\nu^n)_\mathfrak{p}}} + \alpha(\xi_\mathfrak{p}) : G_{(N_\nu^n)_\mathfrak{p}} \longrightarrow \mathcal{E}(n,\nu)$$

is unramified.

In order to complete step 3, it is therefore sufficient to show the existence of an element $\varepsilon \in H^1(G_S, \mathcal{E}(n,\nu)) \subseteq H^1(G_k, \mathcal{E}(n,\nu))$ with $\varepsilon_\mathfrak{p} = \xi_\mathfrak{p}$ for all $\mathfrak{p} \in T$.

The exact sequence

$$H^1(G_S, \mathcal{E}(n,\nu)) \longrightarrow \prod_T H^1(k_\mathfrak{p}, \mathcal{E}(n,\nu)) \xrightarrow{\pi_n} \operatorname{coker}(k_S, T, \mathcal{E}(n,\nu))$$

shows that the obstruction to the existence of such an $\varepsilon$ is the vanishing of $\pi_n(\xi_n)$ with

$$\xi_n = \prod_{\mathfrak{p} \in T_n} \xi_\mathfrak{p} \in \prod_{\mathfrak{p} \in T_n} H^1(k_\mathfrak{p}, \mathcal{E}(n,\nu)).$$

(In the following we denote the sets $T^i$ and $T$ by $T_n^i$ and $T_n$, respectively, in order to indicate on which level the embedding problem is considered.) By lemma 10 we have a canonical injection

(1) $\qquad \operatorname{coker}(k_{S_n}, T_n, \mathcal{E}(n,\nu)) \hookrightarrow \text{III}^1(k_{S_n}, S_n \smallsetminus T_n, \mathcal{E}(n,\nu)')^*,$

where $\mathcal{E}(n,\nu)' = \operatorname{Hom}(\mathcal{E}(n,\nu), \mu_p)$.

Recall that $H^{-2} = H_1$ and that by the induction hypothesis, the solution $\varphi_{n,\nu}$ is proper. Thus we obtain

$$\begin{aligned} H_1(\mathcal{F}(n)/\nu \rtimes G, \mathcal{E}(n,\nu)(-1)) &\cong H^1(\mathcal{F}(n)/\nu \rtimes G, \mathcal{E}(n,\nu)^*(1))^* \\ &\cong H^1(N_\nu^n|k, \mathcal{E}(n,\nu)')^*. \end{aligned}$$

Therefore it exists a canonical isomorphism

(2) $\qquad H^{-2}(\mathcal{F}(n)/\nu \rtimes G, \mathcal{E}(n,\nu)(-1)) \xrightarrow{\sim} H^1(N_\nu^n|k, \mathcal{E}(n,\nu)')^*.$

Now we are going to shrink the obstruction to the existence of a 1-class $\varepsilon$ as above. If $\varphi_{n,\nu}$ is induced by a $\varphi_{m,\nu}$ for $m \geq n$ via a $G$-invariant surjection $\mathcal{F}(m) \twoheadrightarrow \mathcal{F}(n)$, then the inclusion (1) is obviously also true in the form

(1)' $\qquad \operatorname{coker}(k_{S_m}, T_m, \mathcal{E}(n,\nu)) \hookrightarrow \text{III}^1(k_{S_m}, S_m \smallsetminus T_m, \mathcal{E}(n,\nu))^*,$

where $S_m$ and $T_m$ are chosen as above but on the level $N_\nu^m$. Using proposition 7 we choose $m \geq n$ such that an (one) arbitrary chosen element in $H^{-2}(\mathcal{F}(m)/\nu \rtimes G, \mathcal{E}(m,\nu)(-1))$ is annihilated by the map, which is induced by a suitable chosen surjection $\mathcal{F}(m) \twoheadrightarrow \mathcal{F}(n)$. Then we consider the diagrams, in which we write c for coker, $\mathcal{E}_d$ for $\mathcal{E}(d,\nu)$ ($d = m,n$) and $\tilde{S}_n = cs(N_\nu^n|k) \cup T_m$:

$$\begin{array}{ccccccc}
H^1(k_{S_m}|k, \mathcal{E}_m) & \longrightarrow & \prod_{T_m} H^1(k_\mathfrak{p}, \mathcal{E}_m) & \xrightarrow{\pi} & c(k_{S_m}, T_m, \mathcal{E}_m) & \hookrightarrow & \text{III}^1(k_{S_m}, S_m \backslash T_m, \mathcal{E}'_m)^* \\
\downarrow \psi & & \downarrow & & \downarrow & & \downarrow \\
H^1(k_{S_m}|k, \mathcal{E}_n) & \longrightarrow & \prod_{T_m} H^1(k_\mathfrak{p}, \mathcal{E}_n) & \xrightarrow{\pi} & c(k_{S_m}, T_m, \mathcal{E}_n) & \hookrightarrow & \text{III}^1(k_{S_m}, S_m \backslash T_m, \mathcal{E}'_n)^* \\
\downarrow {\scriptstyle inf} & & \| & & \downarrow & & \downarrow \\
H^1(k_{\tilde{S}_n}|k, \mathcal{E}_n) & \longrightarrow & \prod_{T_m} H^1(k_\mathfrak{p}, \mathcal{E}_n) & \xrightarrow{\pi} & c(k_{\tilde{S}_n}, T_m, \mathcal{E}_n) & \hookrightarrow & \text{III}^1(k_{\tilde{S}_n}, \tilde{S}_n \backslash T_m, \mathcal{E}'_n)^*,
\end{array}$$

$$\begin{array}{ccccc}
\mathrm{III}^1(k_{S_m}, S_m\setminus T_m, \mathcal{E}'_m)^* & \xrightarrow{\sim} & H^1(N_\nu^m|k, \mathcal{E}'_m)^* & \xleftarrow{\sim} & H^{-2}(\mathcal{F}(m)/\nu \rtimes G, \mathcal{E}_m(-1)) \\
\downarrow & & \downarrow & & \Big\downarrow \\
\mathrm{III}^1(k_{S_m}, S_m\setminus T_m, \mathcal{E}'_n)^* & \longrightarrow\!\!\!\!\!\rightarrow & H^1(N_\nu^n|k, \mathcal{E}'_n)^* & & \\
\downarrow & & \| & & \Big\downarrow \\
\mathrm{III}^1(k_{\tilde{S}_n}, \tilde{S}_n\setminus T_m, \mathcal{E}'_n)^* & \xrightarrow{\sim} & H^1(N_\nu^n|k, \mathcal{E}'_n)^* & \xleftarrow{\sim} & H^{-2}(\mathcal{F}(n)/\nu \rtimes G, \mathcal{E}_n(-1)).
\end{array}$$

The existence of all maps and the fact that the diagrams are commutative follows from the arguments above and from lemma 12.
Now let
$$\xi_m = \prod_{\mathfrak{p}\in T_m} \xi_\mathfrak{p} \in \prod_{\mathfrak{p}\in T_m} H^1(k_\mathfrak{p}, \mathcal{E}(n,\nu)).$$
be arbitrary. By theorem 7(i) that we can choose a $G$-invariant surjection $\psi: \mathcal{F}(m) \twoheadrightarrow \mathcal{F}(n)$ such that $\pi_n \psi_*(\xi_m) = 0$.

Observe that we did not get precisely what we wanted, because $\varepsilon$ has the required property with respect to the sets of primes $T_m^i$ and $\tilde{S}_n$. Nevertheless, one easily verifies that, if we modify the solution that we have obtained after the shrinking by the cocycle $\varepsilon$ (which now exists), then we obtain a solution satisfying all required properties. This finishes step 3.

**Fourth Step:** *After changing $\varphi_{n,\nu}$ again, there exists a proper solution $\varphi_{n,\nu+1}$ of $(*)$ which satisfies properties (i) and (ii).*

The solution $\varphi_{n,\nu+1}$, which we have obtained in step 3, has almost all properties we need, except that for $\mathfrak{p} \in \mathrm{Ram}(N_{\nu+1}^n|K)\setminus\mathrm{Ram}(N_\nu^n|K)$ the local extension $(N_{\nu+1}^n)_\mathfrak{p}|k_\mathfrak{p}$ might not be (cyclic) totally ramified. But we know that for such a prime $\mathfrak{p}$ the extension $(N_\nu^n)_\mathfrak{p}|k_\mathfrak{p}$ is trivial. In order to get a totally ramified cyclic extension, we have to remove the unramified part of the extension $(N_{\nu+1}^n)_\mathfrak{p}|(N_\nu^n)_\mathfrak{p}$ and to make sure that at places, where new ramification occurs by this procedure, we have (cyclic) totally ramified local extensions.

In order to save the properness of the solution, obtained in step 3, we choose a finite set of primes $T^0 \subseteq cs(N_\nu^n|k) \smallsetminus (\mathrm{Ram}(N_{\nu+1}^n|k) \cup S_p \cup S_\infty)$ such that their decomposition groups $G_\mathfrak{p}(N_{\nu+1}^n|N_\nu^n)$, $\mathfrak{p} \in T^0$ generate $G(N_{\nu+1}^n|N_\nu^n)$.

We want to alter the solution found in step 3 once again using a class $x$ in $H^1(k_S|K, \mathcal{E}(n,\nu))$, where

- $S = cs(N_\nu^n|k) \cup T$ and $T = \mathrm{Ram}(N_{\nu+1}^n|k) \cup S_p \cup S_\infty \cup T^0$,
- for $\mathfrak{p} \in \mathrm{Ram}(N_\nu^n|k) \cup S_p \cup S_\infty \cup T^0$ we have $x_\mathfrak{p} = 0$,
- if the prolongations of $\mathfrak{p}$ to $K$ are in $\mathrm{Ram}(N_{\nu+1}^n|K)\setminus\mathrm{Ram}(N_\nu^n|K)$, then $x_\mathfrak{p} \in H^1_{nr}(k_\mathfrak{p}, \mathcal{E}(n,\nu)) = H^1_{nr}((N_\nu^n)_\mathfrak{p}, \mathcal{E}(n,\nu))$ has the property that

$$\varphi_{n,\nu+1}|_{G_{(N_\nu^n)_\mathfrak{p}}} + x_\mathfrak{p} \in H^1((N_\nu^n)_\mathfrak{p}, \mathcal{E}(n,\nu))$$

is cyclic.

· $x_\mathfrak{p}$ is cyclic for all $\mathfrak{p} \notin T$.

For every $\mathfrak{p} \in S(k)$ such that prolongations of $\mathfrak{p}$ to $K$ are in $\mathrm{Ram}(N_{\nu+1}^n|K)$ but not in $\mathrm{Ram}(N_\nu^n|K)$ we fix one ($\mathfrak{p}$ splits completely in $K|k$) prolongation $\mathfrak{p}_0 \in S(K)$ of $\mathfrak{p}$ to $K$. Let $\eta \in \prod_T H^1(K_\mathfrak{P}, \mathcal{E}(n,\nu))$ be such that

· $\eta_\mathfrak{P} = 0$ if $\mathfrak{P} \cap k \in \mathrm{Ram}(N_\nu^n|k) \cup S_p \cup S_\infty \cup T^0$,
· if $\mathfrak{P} \in \mathrm{Ram}(N_{\nu+1}^n|K)\backslash\mathrm{Ram}(N_\nu^n|K)$ and $\mathfrak{P} \neq (\mathfrak{P} \cap k)_0$, then $\eta_\mathfrak{P} = 0$,
· if $\mathfrak{P} \in \mathrm{Ram}(N_{\nu+1}^n|K)\backslash\mathrm{Ram}(N_\nu^n|K)$ and $\mathfrak{P} = (\mathfrak{P} \cap k)_0$, then
  $\eta_\mathfrak{p} \in H^1_{nr}(K_\mathfrak{P}, \mathcal{E}(n,\nu)) = H^1_{nr}((N_\nu^n)_\mathfrak{P}, \mathcal{E}(n,\nu))$ has the property that
  $$\varphi_{n,\nu+1}|_{G_{(N_\nu^n)_\mathfrak{P}}} + \eta_\mathfrak{P} \in H^1((N_\nu^n)_\mathfrak{P}, \mathcal{E}(n,\nu))$$
  is cyclic.

Applying theorem 13 in the situation where $\Omega$ is the field $N_\nu^n$ and $A = \mathcal{E}(n,\nu)$, we see that in order to finish the proof, it suffices to construct an element $y$ in $H^1(k_S|K, \mathcal{E}(n,\nu))$ with $y_\mathfrak{P} = \eta_\mathfrak{P}$ for all $\mathfrak{P} \in T$. Indeed, by this procedure we get new ramification only at places which are completely decomposed in $N_\nu^n|k$. Hence their decomposition groups are cyclic (by theorem 13) and contained in the $p$-elementary abelian group $G(N_{\nu+1}^n|N_\nu^n) \cong \mathcal{E}(n,\nu)$. Thus the local extensions associated to these new ramification primes are cyclic of order $p$, in particular, totally ramified. Furthermore, by the choice of $T^0$, the new solution remains proper.

Similar to the situation with the class $\varepsilon$ in step 3, the exact sequence

$$H^1(K_S|K, \mathcal{E}(n,\nu)) \longrightarrow \prod_T H^1(K_\mathfrak{p}, \mathcal{E}(n,\nu)) \xrightarrow{\pi_n} \mathrm{coker}(K_S, T, \mathcal{E}(n,\nu))$$

shows that the obstruction to the existence of such a $y$ is $\pi_n(\eta) = 0$. Now we apply the shrinking procedure as in step 3, but the commutative diagrams which are used there have to be modified as follows: Replace in the first diagram $k$ by $K$ and consider instead of the second the following diagram

$$\begin{array}{ccccc}
\mathrm{III}^1(K_{S_m}, S_m\backslash T_m, \mathcal{E}'_m)^* & \xrightarrow{\sim} & H^1(N_\nu^m|K, \mathcal{E}'_m)^* & \xleftarrow{\sim} & H^{-2}(\mathcal{F}(m)/\nu, \mathcal{E}_m(-1)) \\
\downarrow & & \downarrow & & \downarrow \\
\mathrm{III}^1(K_{S_m}, S_m\backslash T_m, \mathcal{E}'_n)^* & \longrightarrow & H^1(N_\nu^n|K, \mathcal{E}'_n)^* & & \\
\downarrow & & \| & & \downarrow \\
\mathrm{III}^1(K_{\tilde{S}_n}, \tilde{S}_n\backslash T_m, \mathcal{E}'_n)^* & \xrightarrow{\sim} & H^1(N_\nu^n|K, \mathcal{E}'_n)^* & \xleftarrow{\sim} & H^{-2}(\mathcal{F}(n)/\nu, \mathcal{E}_n(-1)).
\end{array}$$

Then we use part (ii) theorem 7 instead of part (i).

Therefore, after a further shrinking, we get a class $y$ with the properties above. Now theorem 13 induces the existence of the desired class $x \in H^1(k_S|k, \mathcal{E}(n,\nu))$.

The new solution $\tilde{\varphi}_{n,\nu+1} = {}^x\varphi_{n,\nu+1}$ fulfills condition (i) and (ii), hence step 4 and the proof in the case $p \neq \mathrm{char}(k)$ of theorem 15 are complete.

The proof in the case $p = \mathrm{char}(k)$ is comparatively easy: Again we proceed by induction on $\nu$, where $n$ is arbitrary. The case $\nu = (1,1)$ is trivial. In the case $\nu = (2,1)$ we get an embedding problem with abelian kernel isomorphic to $\mathbb{F}_p[G]^n$, which is properly solvable by proposition 11. In the next induction steps we do not have to care about the properness of the solutions, because they are automatically proper by the Frattini argument. By [17] II§2 prop.3 we have $cd_p G_k = 1$, hence $G_k$ is $p$-projective by [17] I§3 prop.16 and we can solve the embedding problems in all induction steps. Therefore the proof of theorem 15 and also that of theorem 14 is complete. □

In order to deduce the theorem of Šafarevič, we use an argument which goes back to O. ORE [9]. We need two facts from group theory and we recall the following notations:
Let $G$ be a finite non-trivial group, then

$\Phi(G)$ is the intersection of all maximal subgroups of $G$ and is called *Frattini group* of $G$,

$F(G)$ is the composite of all nilpotent normal subgroups of $G$ and is called *Fitting group* of $G$.

The group $\Phi(G)$ is a characteristic subgroup of $G$ and is contained in $F(G)$. The group $F(G)$ is obviously a normal subgroup of $G$. We cite the following two facts, see [4] Kap. III Satz 3.2 (b) and Satz 4.2 (c):

**Proposition 16** *Let $N$ be a normal subgroup of the finite group $G$ such that $N \nsubseteq \Phi(G)$. Then there exists a proper supplement $U$ of $N$ in $G$, i.e. $U \neq G$ and $G = N \cdot U$.*

**Proposition 17** *Let $G$ be a non-trivial finite solvable group. Then $\Phi(G)$ is a proper subgroup of $F(G)$.*

**Proof of Šafarevič's theorem:** Let $F(G)$ be the Fitting subgroup of $G \neq \{1\}$. By the two propositions above, $F(G)$ has a proper (solvable) supplement $U \subsetneq G$, hence there exists a surjection

$$F(G) \rtimes U \twoheadrightarrow G \ .$$

Assuming inductively (on the order of $G$) that $U$ is the Galois group of a finite normal extension of $k$, we obtain the result using theorem 14. □

**Acknowledgment:** The authors want to thank *B.H. MATZAT* for pointing out an error in an earlier version of this article.

# References

[1] Artin, E., Tate, J. *Class Field Theory.* Benjamin New York, Amsterdam 1967

[2] Fried, M.D., Jarden, M. *Field Arithmetic.* Springer 1996

[3] Hoechsmann, K. *Zum Einbettungsproblem.* J. reine u. angew. Math. **229** (1968) 81-106

[4] Huppert, B. *Endliche Gruppen I.* Springer 1967

[5] Išhanov, V.V., Lur'e, B.B., Faddeev, D.K. *The Embedding Problem in Galois Theory.* Trans. of Math. Monographs **165** AMS Providence 1997

[6] Neukirch, J. *Über das Einbettungsproplem der algebraischen Zahlentheorie* Invent. Math. **21** (1973) 59-116

[7] Neukirch, J. *On solvable number fields.* Invent. Math. **53** (1979) 135-164

[8] Neukirch, J. *Algebraische Zahlentheorie.* Springer 1992, Engl.tranl: Algebraic Number Theory. Springer 1998

[9] Ore, O. *Contributions to the theory of groups of finite order.* Duke Math. J. **5** (1939) 431-460

[10] Reichardt, H. *Konstruktion von Zahlkörpern mit gegebener Galoisgruppe von Primzahlpotenzordnung.* J. reine u. angew. Math. **177** (1937) 1-5

[11] Šafarevič, I.R. *On the construction of fields with a given Galois group of order $\ell^\alpha$.* Izv. Akad. Nauk SSSR **18** (1954) 261-296 (in Russian), English translation in Amer. Math. Soc. Transl. **4** (1956) 107-142

[12] Šafarevič, I.R. *On an existence theorem in the theory of algebraic numbers.* Izv. Akad. Nauk SSSR **18** (1954) 327-334 (in Russian), English translation in Amer. Math. Soc. Transl. **4** (1956) 143-150

[13] Šafarevič, I.R. *On the problem of imbedding fields.* Izv. Akad. Nauk SSSR **18** (1954) 389-418 (in Russian), English translation in Amer. Math. Soc. Transl. **4** (1956) 151-183

[14] Šafarevič, I.R. *Construction of fields of algebraic numbers with given solvable groups.* Izv. Akad. Nauk SSSR **18** (1954) 525-578 (in Russian), English translation in Amer. Math. Soc. Transl. **4** (1956) 185-237

[15] Serre, J.-P. *A Course in Arithmetic.* Springer 1973

[16] Serre, J.-P. *Topics in Galois Theory.* Jones and Bartlett Publ. Boston 1992

[17] Serre, J.-P. *Cohomologie Galoisienne.* Lecture Notes in Mathematics **5**, Springer 1964 (Cinquième édition 1994).

[18] Witt, E. *Treue Darstellung Liescher Ringe.* J. reine u. angew. Math. **177** (1937) 152-160

Mathematisches Institut
der Universität Heidelberg
Im Neuenheimer Feld 288
69120 Heidelberg
Germany

Mathematisches Institut
der Universität Heidelberg
Im Neuenheimer Feld 288
69120 Heidelberg
Germany

e-mail: schmidt@mathi.uni-heidelberg.de
        wingberg@mathi.uni-heidelberg.de